\documentclass[12pt]{amsart}

\usepackage{etex}
\usepackage{amscd}
\usepackage{amsmath}
\usepackage{amsfonts}
\usepackage{amssymb}
\usepackage{graphics}
\usepackage{epsfig}
\usepackage{pictex} 

\font\Bbb=msbm10 scaled \magstep 2

\def\P{\hbox{\Bbb P}}

\def\C{\hbox{\Bbb C}}

\def\N{\hbox{\Bbb N}}

\font\midBbb=msbm8

\def\midC{\hbox{\midBbb C}}

\newtheorem{theorem}{\bf Theorem}
\newtheorem{lemma}[theorem]{\bf Lemma}
\newtheorem{remark}[theorem]{\bf Remark}
\newtheorem{example}[theorem]{\bf Example}

\newtheorem{corollary}[theorem]{\bf Corollary}
\newtheorem{definition}[theorem]{\bf Definition}

\numberwithin{equation}{section}
\numberwithin{theorem}{section}

\title{Parameterizing and inverting analytic mappings with~unit Jacobian}

\author{Timur Sadykov}
\address{Plekhanov Russian University of Economics
\newline \indent 115054, Moscow, Russia}
\email{Sadykov.TM@rea.ru}
\thanks{This research was supported by the Russian Science Foundation, project no.~22-21-00556, URL: https://rscf.ru/en/project/22-21-00556/}

\begin{document}

\begin{abstract}
Let $x=(x_1,\ldots,x_n)\in\midC^n$ be a vector of complex variables, denote by $A=(a_{jk})$ a square matrix of size $n\geq 2,$
and let $\varphi\in\mathcal{O}(\Omega)$ be an analytic function defined in a nonempty domain $\Omega\subset\midC.$
We~investigate the family of mappings
$$
f=(f_1,\ldots,f_n):\midC^n\rightarrow\midC^n, \quad f[A,\varphi](x):=x+\varphi(Ax)
$$
with the coordinates
$$
f_j : x \mapsto x_j + \varphi\left(\sum\limits_{k=1}^n a_{jk}x_k\right), \quad j=1,\ldots,n
$$
whose Jacobian is identically equal to a nonzero constant for any~$x$ such that all of~$f_j$ are well-defined.

Let~$U$ be a square matrix such that the Jacobian of the mapping $f[U,\varphi](x)$ is a nonzero constant for any~$x$
and moreover for any analytic function $\varphi\in\mathcal{O}(\Omega).$
We~show that any such matrix~$U$ is uniquely defined, up to a suitable permutation similarity of matrices,
by a partition of the dimension~$n$ into a sum of~$m$ positive integers together with a permutation on~$m$ elements.

For any $d=2,3,\ldots$ we construct $n$-parametric family of square matrices $H(s), s\in\midC^n$
such that for any matrix~$U$ as above the mapping $x+\left((U\odot H(s))x\right)^d$ defined by the Hadamard product $U\odot H(s)$ has unit Jacobian.
We~prove any such mapping to be polynomially invertible and provide an explicit recursive formula for its inverse.

\end{abstract}


\maketitle

\section{Introduction
\label{sec:introduction}}

Let $f=(f_1,\ldots,f_n)$ be an $n$-tuple of analytic functions in $n$~complex variables $x=(x_1,\ldots,x_n)\in\C^n$
defined in a nonempty domain $D\subset\C^n$.
We~will say that~$f$ is a {\it Jacobian mapping} if the determinant of its Jacobian matrix is a nonzero complex number:
\begin{equation}
\label{J}
J(f;x)=
J(f_1,\ldots,f_n; x_1,\ldots,x_n):=
\left|
\begin{array}{ccc}
\frac{\partial f_1}{\partial x_1} & \ldots & \frac{\partial f_n}{\partial x_1} \\
\ldots & \ldots & \ldots \\
\frac{\partial f_1}{\partial x_n} & \ldots & \frac{\partial f_n}{\partial x_n} \\
\end{array}
\right| \in\C^{*} = \C\setminus \{0\}.
\end{equation}
We~will call this determinant {\it the Jacobian of~$f$} as long as this does not lead to ambiguity.

By~the inverse mapping to~$f$ we mean the analytic mapping $f^{-1} : \C^n \rightarrow\C^n$ such that $f\circ f^{-1} = f^{-1}\circ f = {\rm Id}$
in a nonempty domain in~$\C^n.$ The domain where these equalities are valid is in general heavily dependent on the mapping~$f$
and on the domain~$D.$ In~the present paper we do not investigate this dependence since the mappings under study are either
defined by entire functions or are analytically extendible into the whole of~$\C^n$ except for certain singular hypersurface~$\mathcal{H}.$
Such a mapping is uniquely defined by any of its germs which can be further analytically extended into $\C^n\setminus \mathcal{H}.$

If~$n=1,$ then any mapping satisfying~(\ref{J}) is of the form $ax+b$ with $a\in\C^{*}.$
The inverse mapping is also affine linear. From now on, we will disregard this trivial case and assume that $n\geq 2.$

The famous Jacobian Conjecture~\cite{Keller1939} states that a mapping $f=(f_1,\ldots,f_n)$
defined by polynomials $f_j\in\C[x_1,\ldots,x_n]$ is Jacobian
if and only if its inverse is a polynomial mapping, too.
In~any dimension $n\geq 2$ the family of all polynomial mappings $f:\C^n\rightarrow\C^n$ with unit Jacobian
is extremely vast and, judging by known results, has a highly complex structure (see \cite[Chapters~3,5]{Essen} and references therein).
The seminal theorem by Dru$\dot{\rm{z}}$kowski~\cite{Druzkowski1983} ensures that in order to verify the Jacobian Conjecture
it suffices to check its veracity for the very special class of cubic polynomial mappings of the form
$$
\{
(f_1,\ldots,f_n):\C^n\rightarrow\C^n, \,
f_j(x) = x_j + \left( a_{j1}x_1 + \ldots + a_{jn}x_n \right)^3, \, j=1,\ldots,n, \, a_{jk}\in\C
\}.
$$

In~the present paper we investigate a wider family of analytic mappings defined by a square matrix $A=(a_{jk})$ of size $n\geq 2$
together with a function $\varphi\in\mathcal{O}(\Omega)$ that is analytic in a nonempty domain $\Omega\subset\C.$
This family comprises mappings of the form
$$
f=(f_1,\ldots,f_n):\C^n\rightarrow\C^n, \quad f[A,\varphi](x):=x+\varphi(Ax)
$$
with the coordinates
$$
f_j : x \mapsto x_j + \varphi\left(\sum\limits_{k=1}^n a_{jk}x_k\right), \quad j=1,\ldots,n
$$
whose Jacobian is identically equal to a nonzero constant for any~$x$ such that all of~$f_j$ are well-defined.

Let~$U$ be a square matrix such that the Jacobian of the mapping $f[U,\varphi](x)$ is a nonzero constant for any~$x$
and moreover for any analytic function $\varphi\in\mathcal{O}(\Omega).$ We~call such matrices {\it universal}
and show that any universal matrix is uniquely defined, up to a suitable permutation similarity of matrices,
by a partition of the dimension~$n$ into a sum of~$m$ positive integers together with a permutation on~$m$ elements
(see Theorem~\ref{thm:UnivMatrixIsDefinedByPartitionAndPermutation}).

For any $d=2,3,\ldots$ we construct $n$-parametric family of square matrices $H(s), s\in\C^n$
such that for any universal matrix~$U$ the mapping $x+\left((U\odot H(s))x\right)^d$ defined by the Hadamard product $U\odot H(s)$ has unit Jacobian.
We~prove any such mapping to be polynomially invertible and provide an explicit recursive formula for its inverse
(see Theorem~\ref{thm:HadamardProductsOfUniversalAndHomogeneitiesAreGood}).
Being able to explicitly invert a polynomial mapping allows one to compute complex amoebas
of multivariate polynomials with high precision and to investigate their geometric properties~\cite{Bogdanov-Sadykov}.

All of the polynomial mappings constructed in the paper are polynomially invertible.
Moreover, for a universal matrix~$U$ and any analytic function~$\varphi(\cdot)$ the inverse to the mapping $f[U,\varphi]$
is a finite superposition of~$\varphi(\cdot)$ and elementary arithmetic operations acting on the variables $x=(x_1,\ldots,x_n)\in\C^n$
(see Corollary~\ref{cor:inverseOfJMapDefinedByUMatrixIsComposition}).
However, this does not hold for arbitrary matrix~$A$ and arbitrary analytic function~$\varphi(\cdot)$
such that the Jacobian of the mapping $x+\varphi(Ax)$ identically equals~1.
We employ the notion of circulant matrices to construct examples of Jacobian mappings of the form $f[A,\log](x):=x+\log(Ax)$
and show that for a certain choice of a circulant matrix~$A$ the inverse mapping ${f[A,\log]}^{-1}$ is not a finite superposition
of the logarithmic function and elementary arithmetic operations (see Example~\ref{ex:JEquations,LOG,dim=3}).


\section{Notation and Preliminaries
\label{sec:preliminaries}}

Throughout the paper we denote by~$\# I$ the cardinality of a set~$I.$
The scalar product of two vectors $x,y\in\C^n$ is denoted by $\langle x,y \rangle:=x_1 y_1 + \ldots + x_n y_n.$
The transpose of a matrix~$A$ is denoted~$A^T,$ by~$\delta_{jk}$ we denote the Kronecker delta.
For a~square matrix~$A=(a_{jk})$ of size~$n$ and $I\subset \{1,\ldots,n\}$
we denote by~$[A]_I$ the principal minor of~$A$ corresponding to the multi-index~$I,$
that is, the determinant of the square submatrix $(a_{jk})_{j,k\in I}.$
The rows of a matrix~$A$ are referred to as $A_1,\ldots,A_n.$
By~a strictly upper (respectively lower) triangular matrix we mean a matrix with zeros in the main diagonal and below (respectively over)~it.

For matrices $A=(a_{jk})$ and $B=(b_{jk})$ of equal size we denote by~$A\odot B$ the Hadamard (termwise) product
of~$A$ and~$B,$ i.e., the matrix~$(a_{jk}b_{jk}).$ For $d\in\N$ we denote by~$A^{\odot d}$
the $d$th Hadamard power of~$A,$ i.e., the matrix $A\odot\ldots\odot A$ ($d$~copies of~$A$).
The Kronecker product of~$A$ and~$B$ is denoted by~$A\otimes B.$

For $v=(v_1,\ldots,v_n)\in\C^n$ we denote by~$C(v)$ the circulant matrix whose first row is~$v,$ i.e., the matrix
$$
C(v):=
\left(
\begin{array}{ccccc}
v_1   & v_2  & v_3  & \ldots & v_n     \\
v_n   & v_1  & v_2  & \ldots & v_{n-1} \\
\ldots&\ldots&\ldots&\ldots&\ldots \\
v_2   & v_3  & v_4  & \ldots & v_1 \\
\end{array}
\right).
$$
We~denote by~$V(v)$ the Vandermonde-type square matrix
$$
V(v):=
\left(
\begin{array}{ccccc}
v_{1}  & v_{1}^2 & \ldots & v_{1}^{n} \\
v_{2}  & v_{2}^2 & \ldots & v_{2}^{n} \\
\ldots & \ldots  & \ldots & \ldots    \\
v_{n}  & v_{n}^2 & \ldots & v_{n}^{n} \\
\end{array}
\right).
$$
Here we shift the exponents in the standard definition of the Vandermonde
matrix by~$1$ to cover the case when $v_j=0$ which in our setup must yield a zero row in the matrix~$V(v).$

By~a partition of a positive integer~$n$ we mean any vector $(p_1,\ldots,p_m)$ with positive integer entries such that $p_1 + \ldots + p_m = n.$
We~assume that the elements of a partition of an integer are ordered by ascension, i.e., that $p_j \leq p_{j+1}$
unless explicitly stated otherwise. To~simplify the notation, we will often write a partition of an integer in the form "$p_1 + \ldots + p_m$".
We~furthermore will denote the sum of the elements of a partition~$p$ by $|p|:=p_1 + \ldots + p_m.$
By~a partition of a set we mean any representation of this set as a disjoint union of its nonempty subsets.

The permutation of the set of integers $\{1,2,\ldots,m\}$ which takes~$j$ to~$\alpha_j,$ $j=1,\ldots,m$
is denoted by $(\alpha_{1}\alpha_{2}\ldots\alpha_{m}).$ By~$\Pi(\alpha_{1}\alpha_{2}\ldots\alpha_{m})$ we denote
the corresponding permutation matrix, i.e., the matrix such that left multiplication of
$\Pi(\alpha_{1}\alpha_{2}\ldots\alpha_{m})$ with any matrix~$A$ results in permuting the rows of~$A.$
By~$\varepsilon_m$ we denote the unit in the group of permutations of $\{1,2,\ldots,m\},$ i.e., the identity permutation
$\varepsilon_{m}:=(12\ldots m).$

We~denote by $\mathcal{O}(\Omega)$ the ring of functions that are analytic in a nonempty domain $\Omega\subset\C.$
Throughout the paper we define the action of a univariate analytic function~$\varphi(\cdot)\in\mathcal{O}(\Omega)$
on a complex vector $\xi=(\xi_1,\ldots,\xi_n)\in\Omega\times\ldots\times\Omega\subset\C^n$ to be termwise:
$\varphi(\xi_1,\ldots,\xi_n):=(\varphi(\xi_1),\ldots,\varphi(\xi_n)),$ unless otherwise is explicitly stated.
We~denote by~$\P^n$ the complex projective space of dimension~$n.$


\section{Background from linear algebra}
\label{sec:BackgroundFromLinearAlgebra}

In~this section we collect basic facts from linear algebra that will be used in the sequel.

\subsection{Permutation similarity of square matrices.} Throughout the paper the following definition will be used.
\begin{definition}
\label{def:permutationEquivalentMatrices}\rm
Two square matrices~$A$ and~$B$ of equal size are said to be {\it permutation-similar}
if there exists a permutation matrix~$\Pi$ such that $\Pi A \Pi^T = B.$ Here~$\Pi^T$ is the transpose of~$\Pi.$
\end{definition}
In~other words, two matrices are permutation-similar if one can be obtained from the other by first
permuting its rows and then permuting its columns in accordance with one given permutation.
Since the transpose of a permutation matrix is also its inverse, permutation similarity is indeed a special instance of
the general notion of similarity of matrices.
For instance, the matrices
$\left(\begin{array}{cc}a&b\\ c&d \end{array}\right)$ and $\left(\begin{array}{cc}d&c\\ b&a \end{array}\right)$
are permutation-similar since
$$
\left(\begin{array}{cc}a&b\\ c&d \end{array}\right) = \left(\begin{array}{cc}0&1\\1&0\end{array}\right)
\left(\begin{array}{cc}d&c\\ b&a \end{array}\right) \left(\begin{array}{cc}0&1\\1&0\end{array}\right).
$$
We~will also say that any of these matrices can be obtained from the other one by means of the permutation similarity
defined by the only nontrivial permutation $(21)$ of the set of integers $\{1,2\}.$
To~simplify the notation we~will often identify a permutation with the corresponding permutation matrix as long as it does not lead to ambiguity.

Clearly, the set of matrices that are permutation-similar to a given generic square matrix of size~$n$ comprises~$n!$ different matrices.
For instance, the six permutation-similar $3\times 3$-matrices with different entries~$a_{jk}$ are given by
$$
{\small
\left(
\begin{array}{ccc}
 a_{11} & a_{12} & a_{13} \\
 a_{21} & a_{22} & a_{23} \\
 a_{31} & a_{32} & a_{33} \\
\end{array}
\right),\left(
\begin{array}{ccc}
 a_{11} & a_{13} & a_{12} \\
 a_{31} & a_{33} & a_{32} \\
 a_{21} & a_{23} & a_{22} \\
\end{array}
\right),\left(
\begin{array}{ccc}
 a_{22} & a_{21} & a_{23} \\
 a_{12} & a_{11} & a_{13} \\
 a_{32} & a_{31} & a_{33} \\
\end{array}
\right),
}
$$
$$
{\small
\left(
\begin{array}{ccc}
 a_{22} & a_{23} & a_{21} \\
 a_{32} & a_{33} & a_{31} \\
 a_{12} & a_{13} & a_{11} \\
\end{array}
\right),
\left(
\begin{array}{ccc}
 a_{33} & a_{31} & a_{32} \\
 a_{13} & a_{11} & a_{12} \\
 a_{23} & a_{21} & a_{22} \\
\end{array}
\right),
\left(
\begin{array}{ccc}
 a_{33} & a_{32} & a_{31} \\
 a_{23} & a_{22} & a_{21} \\
 a_{13} & a_{12} & a_{11} \\
\end{array}
\right).
}
$$
Since permutations form a group, permutation similarity of matrices is an equivalence relation.
By~definition, the (unordered) set of all principal minors of a matrix is invariant under any permutation similarity.

\subsection{Solution to the Vandermonde-type system of linear algebraic equations.}
The next elementary lemma is included for the sake of completeness and convenience of future reference.

\begin{lemma}
\label{lemma:genSolToVandermondeSystem}
Let $I_0,I_1,\ldots,I_m$ be a partition of the set $\{1,\ldots,n\},$ i.e., $I_0 \cup I_1 \cup \ldots \cup I_m = \{1,\ldots,n\}$
and $I_j \cap I_k = \emptyset$ for $j\neq k.$ Fix $v=(v_1,\ldots,v_n)\in\C^n$ such that

1) $v_j = 0$ if and only if $j\in I_0;$

2) $v_j = v_k$ if and only if $j,k\in I_{\ell}$ for some $\ell\in \{0,\ldots,m\}.$

Let $s=(s_1,\ldots,s_n)\in\C^n.$ The solution space of the system of~$n$ linear equations $(V(v))^{T}s=0$
is the linear subspace of~$\C^n$ defined by the equations
\begin{equation}
\label{kernelOfTransposedVandermondeMatrix}
\sum\limits_{j\in I_k} s_j = 0, \quad k = 1,\ldots, m.
\end{equation}
\end{lemma}
\begin{proof}
Under the assumptions imposed on the vector~$v$ the matrix $V(v)$ has rank~$m.$
Any vector in the linear space~(\ref{kernelOfTransposedVandermondeMatrix}) is in the kernel
of the matrix $(V(v))^{T}$ while the dimension of this space equals~$n-m.$
\end{proof}

The next corollary is an immediate consequence of Lemma~\ref{kernelOfTransposedVandermondeMatrix}
corresponding to the case when $I_0 = \emptyset.$

\begin{corollary}
\label{cor:sumOfSolsToVandermondeSystem=0}
For any $v\in(\C^{*})^n$ the sum of all coordinates of any solution to the system of linear equations with the matrix~$(V(v))^{T}$
equals zero.
\end{corollary}

\subsection{Principal minors and strictly upper triangular matrices.}
Recall that by a strictly upper triangular matrix we mean a matrix with zeros in the main diagonal and below~it.
There is no difference between upper and lower triangular matrices as long as we are concerned with equivalence classes
of permutation-similar matrices.
Indeed, for any upper triangular square matrix~$A$ of size~$n$ the matrix $\Pi(n\ldots 1)A\,\Pi(n\ldots 1)^{T}$ is a lower triangular matrix.
The next lemma gives a necessary and sufficient condition for a matrix to be permutation-similar to a strictly upper triangular matrix.
\begin{lemma}
\label{lemma:zeroPrincipalMinorsEquivToUpperTriangular}
A~square matrix is permutation-similar to a strictly upper triangular matrix if and only if all of its principal minors vanish.
\end{lemma}
\begin{proof}
By~definition, all of the principal minors of a strictly upper triangular matrix vanish.
Since the (unordered) set of all principal minors of a matrix is invariant under a permutation similarity, it follows that
all of the principal minors of a matrix which is permutation-similar to a strictly upper triangular matrix are equal to zero.

We~prove the converse by induction. A~degenerate $2\times 2$ matrix with zeros in its main diagonal
is clearly permutation-similar to a strictly upper triangular matrix.
Assume that the lemma holds for any square matrix of size~$n$ and let~$A$ be a $(n+1)\times(n+1)$ matrix all of whose principal minors vanish.
Denote by~$S$ the $n\times n$ submatrix of~$A$ comprising its first~$n$ rows and first~$n$ columns.
The principal minors of~$S$ are a subset of the set of the principal minors of~$A$ and therefore vanish.
By~the induction hypothesis the matrix~$S$ is permutation-similar to a strictly upper triangular matrix.
Denote by~$\tilde{A}$ the action of this permutation similarity transformation on the matrix~$A.$

Let now~$\tilde{S}$ denote the $n\times n$ submatrix of~$\tilde{A}$ comprising its last~$n$ rows and the last~$n$ columns.
Using the induction hypothesis for the matrix~$\tilde{S}$ as above we conclude that~$\tilde{A}$ is permutation-similar
to a matrix~$\hat{A}$ such that all of nonzero elements of~$\hat{A},$ possibly except for one element, lie over the main diagonal of~$\hat{A}.$
Moreover this only nonzero element below the main diagonal is necessarily in the first column of~$\hat{A}$ but not in the first row.
We~denote this exceptional element by~$a_{k1}$ where~$k$ is the number of the row of the matrix~$\hat{A}$ which contains~it.
Observe that it is perfectly possible that $a_{k1}=0.$
Hence the matrix~$\hat{A}$ is of the following form
$$
\hat{A}=
\left(
\begin{array}{cccccc}
0      & a_{12} & a_{13} & \ldots & a_{1n}  & a_{1\,n+1} \\
0      &  0     & a_{23} & \ldots & a_{2n}  & a_{2\,n+1} \\
\ldots & \ldots & \ldots & \ldots & \ldots  &  \ldots    \\
a_{k1} &   0    &  0     & \ldots & a_{kn}  & a_{k\,n+1} \\
\ldots & \ldots & \ldots & \ldots & \ldots  &  \ldots    \\
0      &   0    &  0     & \ldots &    0    &  a_{n\,n+1}\\
0      &   0    &  0     & \ldots &    0    &      0     \\
\end{array}
\right).
$$
If~$a_{k1}=0$ then the matrix~$\hat{A}$ is strictly upper triangular and permutation-similar to~$A,$ so we are done.
If~$a_{k1}\neq 0,$ denote by~$\breve{A}$ the $k\times k$ submatrix of~$\hat{A}$ comprising its first~$k$ rows and first~$k$ columns.
Arguing as above we conclude that all of the principal minors of~$\breve{A}$ vanish. Induction on~$k$ shows that the matrix~$\breve{A}$
contains at least one zero column 
and hence so does the matrix~$\hat{A}.$ Let~$\ell$ be the number of this zero column.
Acting, if necessary, by the permutation similarity $(1\ell)$ on the matrix~$\hat{A}$ we may without loss of generality
assume that $\ell=1.$ Using the induction hypothesis on the matrix comprising the last~$n$ rows and the last~$n$ columns of~$\hat{A}$ we arrive
at the strictly upper triangular matrix that is permutation-similar to~$A.$
\end{proof}


\section{Good pairs and universal matrices}
\label{sec:goodPairsAndUnivMatrices}

The set of all Jacobian mappings in any dimension that exceeds one is very vast and probably has a highly complex structure.
In~the present paper we only study Jacobian mappings of a very particular form, namely, the mappings $x + \varphi(Ax),$
where $x\in\C^n,$ $A$~is a square $n\times n$-matrix, and $\varphi(\cdot)$ is an analytic function whose action on a complex
vector is defined to be termwise. The interest in such mappings stems from the groundbreaking theorem
by Dru$\dot{\rm{z}}$kowski~\cite{Druzkowski1983} stating that it suffices to investigate the mappings $x + \varphi(Ax)$ defined by
the cubic function $\varphi(\zeta)=\zeta^3$ and all matrices in all dimensions in order to find out whether the Jacobian Conjecture is true or not.
We~would like to study different instances of the function~$\varphi(\cdot)$ and the corresponding matrices which together define Jacobian mappings.
Throughout the paper the following definition is adopted.

\begin{definition}
\label{def:GoodPair}\rm
A~square matrix~$A=(a_{jk})$ of size~$n$ together with a function $\varphi(\zeta)\in\mathcal{O}(\Omega)$ which is
analytic in a nonempty domain $\Omega\subset\C$ are said to form a {\it good pair}
if the Jacobian of the mapping $x\mapsto x+\varphi(Ax)$ with the coordinates
\begin{equation}
\label{generalMapWithFcnPhi}
x_j + \varphi\left(\sum_{k=1}^n a_{jk}x_k\right), \quad j=1,\ldots,n
\end{equation}
is identically equal to a nonzero constant in the domain $\{x\in\C^n : \sum\limits_{k=1}^n a_{jk}x_k \in \Omega, j=1,\ldots,n \}$
which is also assumed to be nonempty.
\end{definition}

By~the above definition the zero matrix forms a good pair with any function that is analytic at the origin.
The identity matrix forms a good pair with the linear function $\varphi(\zeta)=c\cdot\zeta$ for any $c\neq -1.$
A~constant function forms a good pair with any matrix.
In~the sequel we disregard these and similar trivial cases and focus on finding matrices that form good pairs with
"interesting" analytic functions, e.g. polynomials, entire, meromorphic functions, etc.

Denote the rows of the matrix $A=(a_{jk})$ by~$A_j,$ $j=1,\ldots,n.$
Induction on~$n$ yields that the Jacobian $J(A,\varphi;x)$ of the mapping~(\ref{generalMapWithFcnPhi}) is given by
\begin{equation}
\label{universalMatrixEquation}
\begin{array}{c}
J(A,\varphi;x):=
J(x+\varphi(Ax);x) =
1 + \sum\limits_{1 \leq \# I \leq n} [A]_{I} \prod\limits_{j\in I} \varphi'(\langle A_{j},x \rangle) \equiv \\
1+a_{11}\varphi'(\langle A_1,x\rangle)+\ldots+a_{nn}\varphi'(\langle A_n,x\rangle) + \\
\left[A\right]_{1,2}\varphi'(\langle A_1,x\rangle)\varphi'(\langle A_2,x\rangle) + \ldots +
\left[A\right]_{n-1,n}\varphi'(\langle A_{n-1},x\rangle)\varphi'(\langle A_{n},x\rangle) + \ldots + \\
\det(A)\varphi'(\langle A_1,x\rangle)\ldots\varphi'(\langle A_n,x\rangle).
\end{array}
\end{equation}
Here $[A]_{I}$ is the principal minor of the matrix~$A$ defined by the multi-index $I\subset\{1,\ldots,n\}.$

The next lemma implies that it suffices to consider any representative of the permutation similarity class of the matrix~$A$
to investigate the Jacobian of the mapping $x + \varphi(Ax).$

\begin{lemma}
\label{lemma:permutationEquivMatrAreGood}
Let a matrix~$A$ form a good pair with a function~$\varphi(\cdot).$
The following hold:
1) Any constant multiple of the matrix~$A$ also forms a good pair with the function~$\varphi(\cdot).$

\noindent
2) For any matrix~$B$ which is permutation-similar to~$A,$ the pair $(B,\varphi)$ is good.
\end{lemma}
\begin{proof}
1)~Denote by $J_{k}(A,\varphi; x):=\sum\limits_{\#I=k} [A]_{I} \prod\limits_{j\in I} \varphi'(\langle A_{j},x \rangle)$ the "homogeneous component"
of order $k>0$ of the Jacobian $J(A,\varphi;x).$
We~disregard the trivial cases of the zero matrix~$A$ or the constant function~$\varphi(\cdot)$ in which the statement of the lemma is clearly true.
Induction on the rank of the matrix~$A$ yields that $J(A,\varphi;x)\equiv 1$ if and only if all of $J_{k}(A,\varphi;x)$ vanish identically in~$x$
for all $k=1,\ldots,n.$
It~remains to observe that $J_{k}(c A,\varphi; x) \equiv c^k J_{k}(A,\varphi; c x)$ and the lemma follows.

2)~Let~$\Pi$ be a permutation matrix such that $\Pi A \Pi^T = B$ and denote $y=\Pi^T x.$
We~have
$$
J(B,\varphi;x) = J(\Pi A \Pi^T,\varphi; x) = J(x+\varphi(\Pi A \Pi^T x);x) =
$$
$$
J(\Pi y + \varphi(\Pi A y); \Pi y) = J(y+\varphi(Ay);y)
$$
and hence the Jacobian of the mapping $x + \varphi(Ax)$ is a nonzero constant if and only if so is the Jacobian of the mapping $x + \varphi(Bx).$
\end{proof}

It~follows from Lemma~\ref{lemma:permutationEquivMatrAreGood}
that a~matrix which is permutation-similar to a strictly upper triangular matrix
(that is, to a matrix with zeros in the main diagonal and below~it) forms a good pair with any function~$\varphi(\cdot)$.
This fact is coherent with Lemma~\ref{lemma:zeroPrincipalMinorsEquivToUpperTriangular}
which yields that all of the principal minors of such a matrix are zero
and hence the right-hand side of~(\ref{universalMatrixEquation}) is identically equal to~1.
We~now aim at describing all matrices which enjoy this property, summarized in the next definition.

\begin{definition}
\label{def:universalMatrix}\rm
Let~$\Omega$ be a nonempty domain in~$\C.$
We~will say that a square matrix~$A=(a_{jk})$ is {\it universal} if the Jacobian of the mapping~(\ref{generalMapWithFcnPhi})
is identically equal to a nonzero constant for any domain~$\Omega$ and any analytic function $\varphi\in\mathcal{O}(\Omega)$
such that~(\ref{generalMapWithFcnPhi}) is well-defined.
In~other words, a matrix is called universal if it forms a good pair with any analytic function.
\end{definition}

\begin{remark}
\label{rem:matrixWithZeroRowCanBeUniversal}\rm
If~the matrix~$A$ contains a zero row and the point $\zeta=0$ is not in the domain of definition of the function~$\varphi(\zeta)$ then
the mapping~(\ref{generalMapWithFcnPhi}) is not well-defined. For the sake of simplicity we will still say that
the matrix~$A$ is universal if the Jacobian of~(\ref{generalMapWithFcnPhi}) is identically equal to a nonzero constant for
generic analytic function~$\varphi(\cdot).$ In~particular, we will adopt the convention that the zero matrix is universal.
\end{remark}

The only universal matrix in one dimension is the $1\times 1$ zero "matrix".
Throughout the rest of the paper we will only speak of universal matrices in dimension two or higher.

There exist several trivial families of universal matrices.
For instance, a~matrix that is permutation-similar to a strictly upper triangular (or, equivalently, lower triangular) matrix is universal.
Any matrix all of whose rows are equal and whose columns sum up to the zero vector is universal.
It~follows directly from the definition that any constant multiple of a universal matrix is also universal,
i.e., the set of such matrices is a double cone.
Theorem~\ref{thm:UniversalityThroughSMatrix} to be proved below implies that in the bivariate case there are no other universal matrices.
\begin{example}
\label{ex:unversalMatricesInDim2}\rm
By the convention in Remark~\ref{rem:matrixWithZeroRowCanBeUniversal} the matrices
$$
\left(
\begin{array}{cc}
0 & a \\
0 & 0
\end{array}
\right), \quad
\left(
\begin{array}{cc}
0 & 0 \\
a & 0
\end{array}
\right), \quad
\left(
\begin{array}{cc}
a & -a \\
a & -a
\end{array}
\right)
$$
are universal for any $a\in\C$. Moreover, there are no other universal $2\times 2$-matrices.

The mapping~(\ref{generalMapWithFcnPhi}) defined by the latter matrix is given by
$$
\begin{array}{c}
f_1(x_1,x_2) = x_1 + \varphi(a x_1 - a x_2), \\
f_2(x_1,x_2) = x_2 + \varphi(a x_1 - a x_2).
\end{array}
$$
Straightforward computation shows that for $x_j(f_1,f_2)$ defined through
$$
\begin{array}{l}
x_1(f_1,f_2) = f_1 - \varphi(a f_1 - a f_2), \\
x_2(f_1,f_2) = f_2 - \varphi(a f_1 - a f_2)
\end{array}
$$
the compositions $f\circ x$ and $x\circ f$ are both equal to the identity mapping as long as all the involved superpositions are well-defined.
The above mappings are therefore each other's inverses for generic analytic function~$\varphi(\cdot).$
\end{example}

In~the sequel we will often make use of the double conic property of the set of universal matrices and only describe the base of the cone.
Since the number of matrices that are permutation-similar to a given matrix rapidly grows with dimension,
we will only give one representative in the equivalence class of permutation-similar matrices
as justified by Lemma~\ref{lemma:permutationEquivMatrAreGood}~2).

\begin{example}\rm
\label{ex:unversalMatricesInDim3}
Any universal $3\times 3$-matrix is permutation-similar to one of the following matrices for a suitable choice of $a,b,c\in\C:$
$$
\left(
\begin{array}{ccc}
 0 & a & b \\
 0 & 0 & c \\
 0 & 0 & 0 \\
\end{array}
\right), \quad
\left(
\begin{array}{ccc}
 a & b & -a-b \\
 a & b & -a-b \\
 a & b & -a-b \\
\end{array}
\right), \quad
\left(
\begin{array}{ccr}
 0 & a &  b \\
 0 & c & -c \\
 0 & c & -c \\
\end{array}
\right), \quad
\left(
\begin{array}{ccc}
 0 & a & -a \\
 b & c & -c \\
 b & c & -c \\
\end{array}
\right).
$$
It~is straightforward to check that for any univariate analytic function~$\varphi(\cdot)$ the Jacobian of the mapping
defined by the latter matrix, i.e., the mapping
$$
\begin{array}{l}
f_1 = x_1 + \varphi(a x_2 - a x_3), \\
f_2 = x_2 + \varphi(b x_1 + c x_2 - c x_3), \\
f_3 = x_3 + \varphi(b x_1 + c x_2 - c x_3) \\
\end{array}
$$
is identically equal to~1. Its~inverse is given by
$$
\begin{array}{l}
x_1 = f_1 - \varphi(a(f_2-f_3)),  \\
x_2 = f_2 - \varphi(b f_1 + c(f_2-f_3) - b\varphi(a(f_2-f_3))), \\
x_3 = f_3 - \varphi(b f_1 + c(f_2-f_3) - b\varphi(a(f_2-f_3)))  \\
\end{array}
$$
and has unit Jacobian for any function~$\varphi(\cdot)$ such that the above compositions are well-defined.
The mappings defined by the other universal matrices in this example are also straightforward to invert.
\end{example}

\begin{remark}
\label{rem:ZeroColumnReduction}
\rm
Linear algebra arguments imply that it suffices to study universal matrices which contain neither zero rows nor zero columns.
Indeed, if a universal $n\times n$ matrix~$A$ contains a zero row then acting, if necessary, by a permutation similarity transformation,
we may without loss of generality assume that its first row~$A_1$ is the zero vector.
Since~$A$ is universal, the Jacobian of the mapping $f_j = x_j + \varphi(\langle A_j,x\rangle),$ $j=1,\ldots,n$
is identically equal to~$1$ for any $x\in\C^n$ and generic analytic function~$\varphi(\cdot).$
Restricting the Jacobian to $x_1=0$ we conclude that the $(n-1)\times(n-1)$ submatrix of~$A$ comprising the elements in its last $n-1$ rows
and last $n-1$ columns is universal since $J(f_1,\ldots,f_n;x_1,\ldots,x_n)|_{x_1=0} \equiv J(f_2,\ldots,f_n;x_2,\ldots,x_n)$.

Similar arguments show that a universal matrix with a zero column necessarily comprises a
universal submatrix of a smaller size.
In~fact, if the $k$th column of a matrix~$A$ only contains zeros then
$J(f_1,\ldots,f_n;x_1,\ldots,x_n)\equiv J(f_1,\ldots [k] \ldots,f_n; x_1,\ldots [k] \ldots ,x_n),$
where~$[k]$ denotes omission of the element with the number~$k.$

Conversely, any universal matrix of size~$n$ can be completed by the zero column with $n+1$ elements and arbitrary row with~$n$
elements to obtain a universal matrix of size~$n+1.$ This does not hold for rows since the matrix
$
\tiny
\left(
\begin{array}{ccr}
 0 & 0 &  0 \\
 b & a & -a \\
 c & a & -a \\
\end{array}
\right)$
is not universal unless $b=c.$
\end{remark}

\subsection{Row sums in the blocks of a matrix.}
Let~$A=(a_{jk})$ be a~square matrix of size~$n$ with the rows $A_1,\ldots,A_n.$
We~denote by~$m$ the number of different rows in~$A,$ $1\leq m \leq n.$
As~long as we are interested in the equivalence classes of permutation-similar matrices,
we may without loss of generality assume that any two equal rows of~$A$ are adjacent
and furthermore that the rows of~$A$ are ordered by the numbers of their occurrences in the set~$\{A_1,\ldots,A_n\}.$
That is, after performing, if necessary, a permutation similarity transformation of the elements of~$A$ we may without loss of generality assume
that there exists a partition $p=(p_1,\ldots,p_m)$ of the dimension~$n$ such that $p_1\leq \ldots \leq p_m$ and

1) the first~$p_1$ rows of~$A$ are equal to each other and differ from any of the other $n-p_1$ rows of~$A;$

2) the next~$p_2$ rows of~$A$ are equal to each other and differ from any of the remaining $n - p_1 - p_2$ rows of~$A;$

$\ldots$

$m$) the last~$p_m$ rows of~$A$ are equal to each other.

\begin{definition}\rm
\label{def:theOrderedFormOfAMatrix}
We~will say that a matrix~$A$ satisfying the above assumptions 1)-$m$) and such that $p_1\leq \ldots \leq p_m$ is given in its {\it ordered form.}
\end{definition}

\begin{example}\rm
\label{ex:numericUniversalMatrix}
The ordered form of the numeric universal matrix
$$
\left(
\begin{array}{rcccrc}
 -3  & 4  & 2  & 6  &  5  & 1 \\
 -15 & 11 & 8  & 12 & -11 & 7 \\
 -3  & 4  & 2  & 6  &  5  & 1 \\
 -19 & 13 & 10 & 0  & -13 & 9 \\
 -15 & 11 & 8  & 12 & -11 & 7 \\
 -3  & 4  & 2  & 6  &  5  & 1 \\
\end{array}
\right)
$$
is given by
\begin{equation}
\label{numericUniversalMatrixOrdered}
M:=
\left(
\begin{array}{ccrccr}
 0  & 13 & -13 & 9 & 10& -19 \\
 12 & 11 & -11 & 7 & 8 & -15 \\
 12 & 11 & -11 & 7 & 8 & -15 \\
 6  & 4  &  5  & 1 & 2 & -3  \\
 6  & 4  &  5  & 1 & 2 & -3  \\
 6  & 4  &  5  & 1 & 2 & -3  \\
\end{array}
\right).
\end{equation}
In~this example $p_1=1, p_2=2, p_3=3.$ We~denote by~$M_j$ the $j$th row of the matrix~$M.$
For $x\in\C^6$ and generic analytic function~$\varphi(\cdot)$ the inverse of the mapping $x+\varphi(Mx)$
(i.e.,~the solution to the system of equations $f_j = x_j + \varphi(\langle M_{j},x \rangle),$ $j=1,\ldots,6$) is given by
$$
\begin{array}{l}
x_1 = f_1 - \varphi(13 f_2-13 f_3+9 f_4+10 f_5-19 f_6), \\
x_j = f_j - \varphi(-12 \varphi(13 f_2-13 f_3+9 f_4+10 f_5-19 f_6) \\
\phantom{x_j = f_j} +12 f_1+11 f_2-11 f_3+7 f_4+8 f_5-15 f_6), \quad j=2,3, \\
x_j = f_j - \varphi(-9 \varphi(-12 \varphi(13 f_2-13 f_3+9 f_4+10 f_5-19 f_6) \\
\phantom{x_j = f_j}  +12 f_1+11 f_2-11 f_3+7 f_4+8 f_5-15 f_6)-6 \varphi(13 f_2-13 f_3+9 f_4+10 f_5-19 f_6) \\
\phantom{x_j = f_j}  +6 f_1+4 f_2+5 f_3+f_4+2 f_5-3 f_6), \quad\quad\quad j=4,5,6. \\
\end{array}
$$
\end{example}

The ordered form of a matrix is not unique as the blocks of same rows which comprise equally many rows
can be interchanged by the action of a suitable permutation similarity.

Assume now that a square matrix~$A$ is given in its ordered form.
With the above notation, we define the square matrix~$\mathcal{S}(A)=\left(\mathcal{S}_{jk}(A)\right)$ of size~$m$ through
\begin{equation}
\label{matrixOfRowSumsInBlocks}
\mathcal{S}_{jk}(A):=
\sum\limits_{\ell = p_1 + p_2 + \ldots + p_{k-1} + 1}^{p_1 + p_2 + \ldots + p_k}
a_{p_1 + p_2 +\ldots + p_j, \ell} \, ,
\quad j,k=1,\ldots,m.
\end{equation}
For instance, for the numeric universal matrix~(\ref{numericUniversalMatrixOrdered}) the matrix~$\mathcal{S}(M)$
is given by the strictly lower triangular matrix
$$
\mathcal{S}(M)=
\left(
\begin{array}{ccc}
 0  & 0 & 0 \\
 12 & 0 & 0 \\
 6  & 9 & 0 \\
\end{array}
\right).
$$

We~now establish a necessary and sufficient condition for a matrix to be universal.

\begin{theorem}
\label{thm:UniversalityThroughSMatrix}
Let~$A$ be a~square matrix given in its ordered form.
The matrix~$A$ is universal if and only if all of the principal minors of the matrix~$\mathcal{S}(A)$ vanish.
\end{theorem}
\begin{proof}
By~convention the statement is trivially true for the zero matrix~$A.$ In~the sequel we disregard this degenerate case.

We~argue by induction on the number~$m$ of different rows in the matrix~$A.$
If~$m=1,$ i.e., if all rows of~$A$ are equal, then the only principal minors of~$A$ that do not necessarily vanish are its diagonal elements.
It~follows from~(\ref{universalMatrixEquation}) that the Jacobian of the mapping $x+\varphi(Ax)$ is given by
$J(A,\varphi;x)=1+({\rm tr}\,A)\,\varphi'(\langle A_1,x\rangle).$
The matrix~$A$ is universal if and only if $J(A,\varphi;x)\equiv 1$ for generic analytic function~$\varphi(\cdot)$
which is equivalent to ${\rm tr}\,A=0.$
Since all of the rows of~$A$ are equal, its trace is also the sum of all elements in its first row and hence by~(\ref{matrixOfRowSumsInBlocks})
the matrix~$\mathcal{S}(A)$ is the $1\times 1$ zero "matrix".

Let us now assume that the statement of the theorem is valid for all $n\times n$ matrices with at most~$m$ different rows, $m<n.$
Suppose that a~matrix~$A$ of size~$n$ has $m+1$ different rows $\tilde{A}_1,\ldots,\tilde{A}_{m+1}.$
Using~(\ref{universalMatrixEquation}) together with the multilinearity of determinants and the induction hypothesis we conclude that
for any multi-index $I\subset\{1,\ldots,m+1\}$ the principal determinant~$[\mathcal{S}(A)]_{I}$ equals the coefficient
by the product $\prod\limits_{j\in I}\varphi'(\langle\tilde{A}_{j},x\rangle)$ in the Jacobian $J(A,\varphi;x).$
Since~$\varphi(\cdot)$ is a generic analytic function, 
it follows from Lemma~\ref{lemma:genSolToVandermondeSystem}
that the matrix~$A$ is universal if and only if all of these coefficients vanish.
\end{proof}

\begin{remark}\rm
\label{rem:matrixSWithoutOrdering}
It~is of course perfectly possible to define the permutation similarity class of the matrix~$\mathcal{S}(A)$
without transforming the matrix~$A$ into its ordered form. Namely, let $I_1,I_2,\ldots,I_m$ be the partition of the set $\{1,\ldots,n\}$
encoding the subsets of equal rows in the matrix~$A.$ That is, $A_j = A_k$ if and only if $j,k\in I_{\ell}$ for some $\ell=1,\ldots,m.$
For any fixed $i_j\in I_j$ denote
$$
\tilde{\mathcal{S}}_{jk}(A):=
\sum\limits_{\ell\in I_k} a_{i_j, \ell} \, ,
\quad j,k=1,\ldots,m.
$$
Since all rows of the matrix~$A$ that are indexed by the elements of any fixed multiindex~$I_j$ are equal,
it follows that $\tilde{\mathcal{S}}_{jk}(A)$ is independent of the choice of $i_j\in I_j.$
By~construction the matrix $\tilde{\mathcal{S}}(A)=\left(\tilde{\mathcal{S}}_{jk}(A)\right)$ is permutation-similar to~$\mathcal{S}(A).$
\end{remark}

\begin{corollary}
\label{cor:sumsOfPrincipalMinorsOfUnivMatr=0}
For a universal matrix, the sum of its principal minors of any fixed order equals zero.
In~particular, the trace and the determinant of a universal matrix vanish.
\end{corollary}
\begin{proof}
Choosing $\varphi(\zeta) = \zeta^{j}/j,$ $j=2,\ldots, k + 1$ in~(\ref{universalMatrixEquation}) we conclude that
the principal minors of order~$k$ of a universal matrix satisfy the system of linear equations with the matrix~$(V(v))^{T},$
where $v = \left(\prod\limits_{j\in I} \langle A_{j},x \rangle \, : \, \# I = k \right).$ Disregarding the trivial case when the matrix~$A$
contains zero rows we may always choose $x\in\C^n$ such that $v\in(\C^{*})^{\binom{n}{k}}.$
Using Corollary~\ref{cor:sumOfSolsToVandermondeSystem=0} we arrive at the conclusion of the lemma.
\end{proof}

We~now introduce notation to be used in the description of a general universal matrix in any dimension.
For $a=(a_1,\ldots,a_{\ell})\in\C^{\ell}$ we will denote by~$S_{k\ell}[a]$ the $k\times\ell$ matrix
with~$k$ equal rows, each row given by $(a_1,a_2,\ldots,a_{\ell}):$
\begin{equation}
\label{blockWithSameRows}
S_{k\ell}[a]:=
\left(
\begin{array}{cccc}
a_1    & a_2    & \ldots & a_{\ell}    \\
\ldots & \ldots & \ldots & \ldots \\
a_1    & a_2    & \ldots & a_{\ell}    \\
\end{array}
\right).
\end{equation}
Any matrix which is a special instance of~(\ref{blockWithSameRows}) will be called {\it a matrix of type~$S.$}
Furthermore, we denote by~$Z_{k\ell}[a]$ the following $k\times\ell$ matrix with~$k$ equal rows whose elements sum up to zero:
\begin{equation}
\label{blockWithSameRowsZeroColmunSum}
Z_{k\ell}[a]:=
\left(
\begin{array}{ccccc}
a_1    & a_2    & \ldots & a_{\ell-1} & -a_1 - a_2 - \ldots - a_{\ell-1} \\
\ldots & \ldots & \ldots & \ldots  & \ldots                        \\
a_1    & a_2    & \ldots & a_{\ell-1} & -a_1 - a_2 - \ldots - a_{\ell-1} \\
\end{array}
\right), \quad {\rm if\ } \ell\geq 2,
\end{equation}
and define~$Z_{k\ell}[a]$ to be the zero column if $\ell=1.$
We~will call any matrix which is a special case of~(\ref{blockWithSameRowsZeroColmunSum}) {\it a matrix of type~$Z.$}
That is, a type~$Z$ matrix is any matrix all of whose rows are equal and whose columns sum up to the zero vector.
We~remark that the matrix~$Z_{k\ell}[a]$ does not depend on the last coordinate of the vector~$a.$

We~now define the $k\times\ell$ matrix $B_{k\ell}^{\sigma}[a]$ whose elements depend on the coordinates of a vector
$a=(a_1,\ldots,a_{\ell})\in\C^{\ell}$ and a binary parameter~$\sigma$ as follows:
$$
B_{k\ell}^{\sigma}[a]:=
\left\{
\begin{array}{r}
Z_{k\ell}[a], \quad \text{if\ } \sigma = 0, \\
S_{k\ell}[a], \quad \text{if\ } \sigma = 1. \\
\end{array}
\right.
$$

For $m\in\N,$ $m\geq 2$ denote by~$\tau_m$ the $m\times m$ strictly lower triangular matrix with units under the diagonal:
$$
\tau_m :=
\left(
\begin{array}{ccccccc}
0 & 0 &    \ldots & 0 & 0 \\
1 & 0 &    \ldots & 0 & 0 \\
\ldots  &\ldots &\ldots &\ldots &\ldots \\
1 & 1 &    \ldots & 0 & 0 \\
1 & 1 &    \ldots & 1 & 0 \\
\end{array}
\right).
$$
We~set~$\tau_1$ to be the $1\times 1$ zero "matrix".

The next theorem gives a complete description of all universal matrices in any dimension.

\begin{theorem}
\label{thm:UnivMatrixIsDefinedByPartitionAndPermutation}
1) Let~$A$ be a universal matrix of size~$n\geq 2.$
There exist an integer partition $p=(p_1,\ldots,p_m)$ of~$n,$ a~permutation matrix~$\Pi$ of size~$m,$ and $a^{(j,k)}\in\C^{p_k},$ $j,k=1,\ldots,m,$
such that the matrix~$A$ is permutation-similar to the block matrix
\begin{equation}
\label{canonicalFormOfUMatrix}
U(p,\Pi):=
\left(
\begin{array}{cccc}
B_{p_1 p_1}^{\eta_{11}}[a^{(1,1)}]  &  B_{p_1 p_2}^{\eta_{12}}[a^{(1,2)}]  & \ldots & B_{p_1 p_m}^{\eta_{1m}}[a^{(1,m)}] \\
B_{p_2 p_1}^{\eta_{21}}[a^{(2,1)}]  &  B_{p_2 p_2}^{\eta_{22}}[a^{(2,2)}]  & \ldots & B_{p_2 p_m}^{\eta_{2m}}[a^{(2,m)}] \\
      \ldots            &        \ldots            & \ldots &          \ldots        \\
B_{p_m p_1}^{\eta_{m1}}[a^{(m,1)}]  &  B_{p_m p_2}^{\eta_{m2}}[a^{(m,2)}]  & \ldots & B_{p_m p_m}^{\eta_{mm}}[a^{(m,m)}] \\
\end{array}
\right),
\end{equation}
where $\eta = (\eta_{jk}),$ $1\leq j,k \leq m$ is the matrix defined by $\eta = \Pi \tau_m \Pi^T.$

2) The universal matrix $U(p,\Pi)$ depends on $m\left(n-\frac{m+1}{2}\right)$ complex parameters.
\end{theorem}
\begin{proof}
By~Theorem~\ref{thm:UniversalityThroughSMatrix} and Lemma~\ref{lemma:zeroPrincipalMinorsEquivToUpperTriangular}
the matrix~$A$ is universal if and only if there exists a permutation matrix~$\Pi$ (depending on~$A$)
such that the product $\Pi^T \mathcal{S}(A) \Pi$ is a strictly lower triangular matrix.
Thus the matrix~$A$ comprises~$m^2$ blocks of sizes~$p_j \times p_k,$ $j,k=1,\ldots,m,$ each block (indexed by $j,k$)
consisting of equal rows and having zero column sum if and only if the element in the $j$th row and the $k$th column of the matrix
$\Pi\tau_m \Pi^T$ equals~$1.$ Denoting these blocks by $B_{p_j p_k}^{\eta_{jk}}[a^{(j,k)}]$ we arrive at~(\ref{canonicalFormOfUMatrix}).

It~remains to compute the number of parameters in the matrix $U(p,\Pi).$
This number is independent of the permutation~$\Pi$ and it is therefore sufficient to consider the case when~$\Pi$ is the identity permutation.
By~definition, a matrix of type~$Z$ and size~$p_j\times p_k$ depends on~$p_k - 1$ algebraically independent complex parameters
while a matrix of type~$S$ and of the same size depends on~$p_k$ such parameters.
Hence the total number of parameters of a universal matrix defined by an integer partition $p=(p_1,\ldots,p_m)$
of the dimension~$n$ equals
$$
\sum\limits_{j=1}^{m} (n-j) = m\left(n-\frac{m+1}{2}\right).
$$
\end{proof}

\begin{corollary}
\label{cor:UnivMatrixIsNilpotent}
Universal matrices are nilpotent.
\end{corollary}
\begin{proof}
The explicit construction of a general universal matrix immediately yields nilpotency of the matrix $U(p,\varepsilon_m)$ defined by any partition
$p=(p_1,\ldots,p_m)$ and the identity permutation of~$m$ elements. Since the property of being nilpotent is invariant under permutation similarity,
the claim holds for any universal matrix.
\end{proof}

Example~\ref{ex:goodMatr,dim=4,deg=2,rank=2,centrosymmetric,idempotent} shows that a matrix that forms a good pair with a particular fixed function
need not in general be nilpotent. However, for certain functions this might be the case, as shown in Corollary~\ref{cor:UMatrices=GoodForExp}
for the exponential function.

\begin{corollary}
\label{cor:KroneckerProductOfUniMatrixAndAnyMatrixIsUni}
For any universal matrix~$U$ and arbitrary square matrix~$M$ the Kronecker products $U\otimes M$ and $M\otimes U$ are universal.
\end{corollary}
\begin{proof}
The explicit construction~(\ref{canonicalFormOfUMatrix}) of a general universal matrix implies that the Kronecker product $M\otimes U$ is universal, too.
Indeed, a constant multiple of a matrix of type~$S$ is again a matrix of type~$S.$ The same holds for a matrix of type~$Z.$
Thus the Kronecker product of a matrix which has the form~(\ref{canonicalFormOfUMatrix})
with arbitrary matrix is again a matrix of the form~(\ref{canonicalFormOfUMatrix}).
Although Kronecker's product is in general not commutative, $A\otimes B$ and $B\otimes A$ are permutation-similar
for any square matrices~$A$ and~$B$ and hence by Lemma~\ref{lemma:permutationEquivMatrAreGood}~2) the matrix $U\otimes M$ is universal.
\end{proof}

\begin{remark}
\label{rem:KroneckerProductOfGoodMatrices}\rm
If~$(A,\zeta^d)$ is a good pair then for any square matrix~$B$ the Kronecker products $A\otimes B$ and $B\otimes A$
form good pairs with the function~$\zeta^d.$ This is not true for arbitrary function~$\varphi(\cdot)$ as it fails, for instance, for
$\varphi(\zeta)=\log\zeta,$ see Example~\ref{ex:JEquations,LOG,dim=3}.
\end{remark}

In~the next example we construct all universal $4\times 4$-matrices that are nontrivial in the sense of Remark~\ref{rem:ZeroColumnReduction}.
Each of these matrices corresponds to a partition of the dimension, i.e., of the integer $n=4$.
The order of the summands in a partition is unimportant since any rearrangement of the summands
corresponds to the action of a suitable permutation similarity on the universal matrix in question.
We~will assume that the summands in a partition of an integer are in ascending order unless otherwise is explicitly stated.
To~avoid triple index notation that is unnecessarily complex in low dimensions, we give the answer in terms of the elements of
a~generic matrix $A=(a_{jk}).$ The notation of~(\ref{canonicalFormOfUMatrix}) is adopted in the 8-dimensional Example~\ref{ex:complexUMatrix}.

\begin{example}
\label{exAllUMatrsInDim4}\rm
Let $A=(a_{jk}),$ $j,k=1,\ldots,4$ be a universal matrix of size~$4$ and let $x=(x_1,x_2,x_3,x_4)\in\C^4.$
Recall that we~denote by $A_j$ the $j$th row of the matrix~$A$ and that $[A]_{j_1,\ldots,j_d}$ stands for the principal
minor of the matrix~$A$ which corresponds to the rows and columns with the indices $j_1,\ldots,j_d.$
For analytic function~$\varphi(\cdot),$ denote $\psi_{j}:=\varphi'(\langle A_{j},x \rangle),$ $j=1,\ldots,4.$
By~(\ref{universalMatrixEquation}) the matrix~$A$ is universal if and only if the equality
\begin{equation}
\label{JOfUMatrInDim=4}
\begin{array}{c}
a_{11}\psi_{1} + a_{22}\psi_{2} + a_{33}\psi_{3} + a_{44}\psi_{4} + \sum\limits_{1\leq j<k\leq 4} \psi_{j}\psi_{k} [A]_{j,k} \\
+ \psi_{1}\psi_{2}\psi_{3} [A]_{1,2,3} + \psi_{1}\psi_{2}\psi_{4} [A]_{1,2,4}
+ \psi_{1}\psi_{3}\psi_{4} [A]_{1,3,4} + \psi_{2}\psi_{3}\psi_{4} [A]_{2,3,4} \\
+ \psi_{1}\psi_{2}\psi_{3}\psi_{4} \det A
=0
\end{array}
\end{equation}
holds identically for a generic univariate analytic function~$\varphi(\cdot)$ and any $x\in\C^4.$

If~the matrix~$A$ contains a zero row then by Remark~\ref{rem:ZeroColumnReduction} it can be reduced to a universal matrix of smaller size.
We~do not address this trivial case in the present example since it can be reduced to the matrices given in Example~\ref{ex:unversalMatricesInDim3}.
Throughout the rest of this example we assume that every row of~$A$ contains at least one nonzero element.
For $k=1,2,3,\ldots$ let $\varphi_k(\zeta):=\frac{\zeta^{k+1}}{k+1}$ and consider~(\ref{JOfUMatrInDim=4}) for $\varphi=\varphi_k.$
Under the above assumptions the homogeneous components of the left-hand side of~(\ref{JOfUMatrInDim=4}) must vanish simultaneously.

{\it The trivial partition 4=4} corresponds to the case when all rows of the matrix~$A$ are equal.
In~this case all of the principal minors of size~$2$ and higher vanish and hence
the left-hand side of~(\ref{JOfUMatrInDim=4}) equals (for $\varphi=\varphi_k$)
$$
(a_{11}+a_{22}+a_{33}+a_{44})(\langle A_1,x\rangle)^k.
$$
So,~$\rm{tr}\, A=0$ and thus the $4\times 4$-matrix all of whose rows are equal is universal if and only if it has the form
$$
U(4,\varepsilon_1):= Z_{44}[(a_{11}, a_{12}, a_{13}, a_{14})]=
\left(
\begin{array}{cccc}
a_{11} & a_{12} & a_{13} & -a_{11}-a_{12}-a_{13} \\
a_{11} & a_{12} & a_{13} & -a_{11}-a_{12}-a_{13} \\
a_{11} & a_{12} & a_{13} & -a_{11}-a_{12}-a_{13} \\
a_{11} & a_{12} & a_{13} & -a_{11}-a_{12}-a_{13} \\
\end{array}
\right).
$$

{\it The partition 4=1+3} encodes a $4\times 4$ universal matrix~$A$ with three equal rows, the remaining row being different from these.
By letting a suitable permutation similarity act, if necessary, on the matrix~$A,$ we may without loss of generality assume that the three
last rows of~$A$ are equal, i.e., $A_2=A_3=A_4.$ Under this assumption any of the principal minors of~$A$ of order~$3$ vanishes and so does
the determinant of~$A.$ The only principal minors of~$A$ of order~$2$ that are not identically equal to zero are $[A]_{1,2},[A]_{1,3},$
and $[A]_{1,4}.$

Choosing in~(\ref{JOfUMatrInDim=4}) the generic analytic function $\varphi=\varphi_k$ as above and using Lemma~\ref{lemma:genSolToVandermondeSystem},
we conclude that
\begin{equation}
\label{system for partition 1+3}
\begin{array}{c}
a_{11}=0, \\
a_{22} + a_{33} + a_{44}=0, \\
{[}A]_{12} + {[}A]_{13} + {[}A]_{14} = 0.
\end{array}
\end{equation}
Since by assumption $a_{21}=a_{31}=a_{41,}$ the first and the last equations in~(\ref{system for partition 1+3}) imply that
$a_{21}(a_{12}+a_{13}+a_{14})=0.$
If~$a_{21}=0$ then all elements of the first column of~$A$ are zero and by Remark~\ref{rem:ZeroColumnReduction}
it can be reduced to a universal matrix of smaller size.
If~$a_{12}+a_{13}+a_{14}=0$ then the matrix~$A$ has the form
$$
U(1+3,\varepsilon_2):=
\left(
\begin{array}{cccc}
Z_{11}[a_{11}] & Z_{13}[(a_{12},a_{13},a_{14})] \\
S_{31}[a_{21}] & Z_{33}[(a_{22},a_{23},a_{24})] \\
\end{array}
\right)=
\left(
\begin{array}{cccc}
   0   & a_{12} & a_{13} & -a_{12}-a_{13} \\
a_{21} & a_{22} & a_{23} & -a_{22}-a_{23} \\
a_{21} & a_{22} & a_{23} & -a_{22}-a_{23} \\
a_{21} & a_{22} & a_{23} & -a_{22}-a_{23} \\
\end{array}
\right).
$$

{\it The partition 4=2+2} stands for a $4\times 4$-matrix $A=(a_{jk})$ such that its first two rows are equal and so are its last two rows.
We~argue as above making use of~(\ref{JOfUMatrInDim=4}) and Lemma~\ref{lemma:genSolToVandermondeSystem}
to conclude that such a matrix is universal if and only if its elements satisfy the following equations:
$$
\left\{
\begin{array}{c}
a_{11} + a_{22}=0, \\
a_{33} + a_{44} = 0, \\
(a_{13} + a_{14})(a_{31} + a_{32}) = 0. \\
\end{array}
\right.
$$
Any solution to this system of equations yields a matrix which is permutation-similar to the universal matrix
$$
U(2+2,\varepsilon_2)=
\left(
\begin{array}{cccc}
Z_{22}[(a_{11},a_{12})] & Z_{22}[(a_{13},a_{14})] \\
S_{22}[(a_{31},a_{32})] & Z_{22}[(a_{33},a_{34})] \\
\end{array}
\right)=
\left(
\begin{array}{crcc}
a_{11} & -a_{11} & a_{13} & -a_{13} \\
a_{11} & -a_{11} & a_{13} & -a_{13} \\
a_{31} &  a_{32} & a_{33} & -a_{33} \\
a_{31} &  a_{32} & a_{33} & -a_{33} \\
\end{array}
\right).
$$

Straightforward computation shows that the mapping defined by the above matrix together with a univariate function~$\varphi(\cdot),$
i.e., the mapping
$$
\begin{array}{l}
f_j = x_j + \varphi(a_{11}x_1 - a_{11}x_2 + a_{13}x_3 - a_{13}x_4), \quad j=1,2,\\
f_j = x_j + \varphi(a_{31}x_1 + a_{32}x_2 + a_{33}x_3 - a_{33}x_4), \quad j=3,4 \\
\end{array}
$$
indeed has unit Jacobian for any $x\in\C^4.$ The inverse of this mapping is given by
$$
{\small
\begin{array}{l}
x_j = f_j-\varphi\left(a_{11} \left(f_1-f_2\right)+a_{13} \left(f_3-f_4\right)\right), \quad j=1,2, \\
x_j = f_j-\varphi\left(a_{31} f_1+a_{32} f_2 -\left(a_{31}+a_{32}\right)
          \varphi\left(a_{11} \left(f_1-f_2\right)+a_{13} \left(f_3-f_4\right)\right)+ a_{33} \left(f_3-f_4\right)\right), \\
         \hfill  j=3,4. \\
\end{array}
}
$$

{\it The partition 4=1+1+2} corresponds to a $4\times 4$-matrix $A=(a_{jk})$ whose last two rows are equal.
We~argue as above making use of~(\ref{JOfUMatrInDim=4}) and Lemma~\ref{lemma:genSolToVandermondeSystem}
to conclude that such a matrix is universal if and only if it is permutation-similar to a matrix satisfying the following equations:
$$
\left\{
\begin{array}{c}
a_{11}=a_{22}=0, \\
a_{33} + a_{44} = 0, \\
a_{31}(a_{13} + a_{14}) = 0, \\
a_{32}(a_{23} + a_{24}) = 0, \\
a_{12}a_{31}(a_{23} + a_{24}) = 0. \\
\end{array}
\right.
$$
The only nontrivial (i.e., corresponding to a matrix with no zero rows or columns)
solution to this system of equations yields the following universal matrix:
$$
U(1+1+2,\varepsilon_3):=
\left(
\begin{array}{cccc}
   0   & a_{12} & a_{13} & -a_{13} \\
   0   &   0    & a_{23} & -a_{23} \\
a_{31} & a_{32} & a_{33} & -a_{33} \\
a_{31} & a_{32} & a_{33} & -a_{33} \\
\end{array}
\right).
$$

{\it Finally, the partition 4=1+1+1+1} corresponds to a universal $4\times 4$-matrix~$A$ all of whose rows are different.
Choosing $\varphi=\varphi_k$ in~(\ref{JOfUMatrInDim=4}) and using Lemma~\ref{lemma:genSolToVandermondeSystem}
for the homogeneous components of the left-hand side of~(\ref{JOfUMatrInDim=4}) we conclude that all
of the principal minors of~$A$ vanish.
By~Lemma~\ref{lemma:zeroPrincipalMinorsEquivToUpperTriangular} such a matrix is permutation-similar
to a generic strictly upper triangular matrix with no relations between its elements.
It~is therefore trivial in the sense of Remark~\ref{rem:ZeroColumnReduction}.
\end{example}

\begin{remark}\rm
\label{rem:powerOfAUniversalMatrixIsUniversal}
Any positive integer power of a universal matrix is universal, too.
However, this does not necessarily hold for a matrix that forms a good pair with a given fixed function.
For instance, the matrix~(\ref{LogSpecificGood3x3-matrix}) forms a good pair with the logarithmic function, but its square does not.
\end{remark}

\begin{example}\rm
\label{ex:complexUMatrix}
Let $n=8$ and define the partition~$p$ of~$n$ to be $p=(1,2,2,3).$ Since~$p$ has~$4$ elements, it follows that $m=4.$
The strictly lower triangular binary matrix with units below the main diagonal is given by
$$
\tau_4 =
\left(
\begin{array}{ccccccc}
0 & 0 & 0 & 0 \\
1 & 0 & 0 & 0 \\
1 & 1 & 0 & 0 \\
1 & 1 & 1 & 0 \\
\end{array}
\right).
$$
We~take a random permutation~$\alpha$ of the set $\{1,2,3,4\}$ to be, say, $\alpha=(1423).$
The action of the permutation similarity defined by~$\alpha$ on the matrix~$\tau_4$ is given by
$$
\Pi(\alpha) \tau_{4} \Pi(\alpha)^{T}=
\left(
\begin{array}{cccc}
 0 & 0 & 0 & 0 \\
 1 & 0 & 0 & 1 \\
 1 & 1 & 0 & 1 \\
 1 & 0 & 0 & 0 \\
\end{array}
\right).
$$
The universal matrix defined by the partition~$p$ and the permutation~$\alpha$ in accordance with~(\ref{canonicalFormOfUMatrix})
is given by
$$
\tiny{
U(1+2+2+3,(1423))=
\left(
\begin{array}{crrrrrrc}
 0             & a_{1}^{(1,2)} & -a_{1}^{(1,2)} & a_{1}^{(1,3)} & -a_{1}^{(1,3)} & a_{1}^{(1,4)} & a_{2}^{(1,4)} & -a_{1}^{(1,4)}-a_{2}^{(1,4)} \\
 a_{1}^{(2,1)} & a_{1}^{(2,2)} & -a_{1}^{(2,2)} & a_{1}^{(2,3)} & -a_{1}^{(2,3)} & a_{1}^{(2,4)} & a_{2}^{(2,4)} & a_{3}^{(2,4)} \\
 a_{1}^{(2,1)} & a_{1}^{(2,2)} & -a_{1}^{(2,2)} & a_{1}^{(2,3)} & -a_{1}^{(2,3)} & a_{1}^{(2,4)} & a_{2}^{(2,4)} & a_{3}^{(2,4)} \\
 a_{1}^{(3,1)} & a_{1}^{(3,2)} &  a_{2}^{(3,2)} & a_{1}^{(3,3)} & -a_{1}^{(3,3)} & a_{1}^{(3,4)} & a_{2}^{(3,4)} & a_{3}^{(3,4)} \\
 a_{1}^{(3,1)} & a_{1}^{(3,2)} &  a_{2}^{(3,2)} & a_{1}^{(3,3)} & -a_{1}^{(3,3)} & a_{1}^{(3,4)} & a_{2}^{(3,4)} & a_{3}^{(3,4)} \\
 a_{1}^{(4,1)} & a_{1}^{(4,2)} & -a_{1}^{(4,2)} & a_{1}^{(4,3)} & -a_{1}^{(4,3)} & a_{1}^{(4,4)} & a_{2}^{(4,4)} & -a_{1}^{(4,4)}-a_{2}^{(4,4)} \\
 a_{1}^{(4,1)} & a_{1}^{(4,2)} & -a_{1}^{(4,2)} & a_{1}^{(4,3)} & -a_{1}^{(4,3)} & a_{1}^{(4,4)} & a_{2}^{(4,4)} & -a_{1}^{(4,4)}-a_{2}^{(4,4)} \\
 a_{1}^{(4,1)} & a_{1}^{(4,2)} & -a_{1}^{(4,2)} & a_{1}^{(4,3)} & -a_{1}^{(4,3)} & a_{1}^{(4,4)} & a_{2}^{(4,4)} & -a_{1}^{(4,4)}-a_{2}^{(4,4)} \\
\end{array}
\right).
}
$$
The above matrix is given in its ordered form.
Straightforward computation shows that for any~$a_{j}^{(k,\ell)}\in\C,$ $k,\ell=1,\ldots,4$ the matrix $U(1+2+2+3,(1423))$ is indeed universal.
By~Lemma~\ref{lemma:permutationEquivMatrAreGood} any matrix that is permutation-similar to it is universal, too.
The explicit inverse mapping of the Jacobian mapping defined by this matrix and generic analytic function~$\varphi(\cdot)$ (i.e., the solution to the system
of equations $f = x + \varphi(U(1+2+2+3,(1423))x)$ with $x,f\in\C^8$) is too cumbersome to display.
Its~last coordinate is given by
$$
{\small
\begin{array}{c}
x_{8} =
f_{8}-\varphi\left(a_{1}^{(4,1)}
\left(f_1-\varphi\left( a_{1}^{(1,2)}(f_2-f_3) +  a_{1}^{(1,3)} (f_4-f_5) +
 a_{1}^{(1,4)} (f_6-f_8) +  a_{2}^{(1,4)}(f_7-f_8) \right)\right) \right. \\
\left.  + a_{1}^{(4,2)} (f_2-f_3) +  a_{1}^{(4,3)} (f_4-f_5) +  a_{1}^{(4,4)} (f_6-f_8) +  a_{2}^{(4,4)} (f_7-f_8) \right).
\end{array}
}
$$
\end{example}

For certain functions~$\varphi(\cdot)$ the only matrices that form good pairs with~$\varphi(\cdot)$ are the universal matrices.
For instance, this is the case for $\varphi(\zeta)=\exp(\zeta).$

\begin{corollary}
\label{cor:UMatrices=GoodForExp}
A~matrix is universal if and only if it forms a good pair with the exponential function~$\exp(\zeta).$
\end{corollary}
\begin{proof}
By~definition, a universal matrix forms a good pair with any analytic function, in particular, with~$\exp(\zeta).$
Conversely, let a~matrix~$A$ form a good pair with the exponential function. Acting, if necessary, by a permutation similarity on~$A,$
and using Lemma~\ref{lemma:permutationEquivMatrAreGood}, we may without loss of generality assume that~$A$ is given in its ordered form.
Due~to the linear independence of the exponential functions with different arguments, vanishing of the Jacobian of the mapping $x+\exp(Ax)$ implies
that all of the principal minors of the matrix~$\mathcal{S}(A)$ equal zero. By~Theorem~\ref{thm:UniversalityThroughSMatrix} the matrix~$A$
is universal.
\end{proof}

Lemma~\ref{lemma:parametrizationOfJMatrices,dim=2} together with Examples~\ref{ex:JEquations,LOG,dim=3},
\ref{ex:goodMatr,dim=4,deg=2,rank=2,centrosymmetric,idempotent}, and~\ref{ex:dim=4,deg=2,rank=2,rational parametrization}
show that the conclusion of Corollary~\ref{cor:UMatrices=GoodForExp} does not hold for an arbitrary function, i.e., that the
set of matrices that form good pairs with a given function is in general much bigger than the set of universal matrices.


\section{The inverse of a Jacobian mapping defined by a universal matrix}

We~now aim at systematically inverting mappings of the form $x+\varphi(Ux)$ with~$U$ being a universal matrix.
It~turns out that the classical Newton iterations are of use here.
For the sake of completeness we recall the following standard definition.

\begin{definition}\rm
\label{def:NewtonInverse}
Let~$A$ be a $n\times n$-matrix, and let~$\varphi(\cdot)$ be a univariate analytic function defined in a nonempty domain in the complex plane.
The sequence of Newton's iterations $\nu_j : \C^n \rightarrow \C^n$ is defined recursively through
\begin{equation}
\label{NewtonIterations}
\nu_0 = {\rm Id}, \quad \nu_{j} = \nu_{j}[A,\varphi](f) := f - \varphi(A \nu_{j-1}(f)), \, j=1,2,3,\ldots
\end{equation}
Here $f=(f_1,\ldots,f_n)\in\C^n,$ $A \nu_{j}(f)$ is the product of the matrix~$A$ and the vector~$\nu_{j}(f),$
while the function~$\varphi(\cdot)$ acts on a vector componentwise, i.e., $\varphi(\xi_1,\ldots,\xi_n) := (\varphi(\xi_1),\ldots,\varphi(\xi_n))$.
\end{definition}

Throughout the rest of the paper we~assume that all of~$\nu_{j}$ are defined in a certain domain in~$\C^n.$
This assumption is satisfied if the function~$\varphi(\cdot)$ is entire, in particular, if $\varphi(\zeta)=\zeta^d,$ $d\in\N.$
However, the iterations~(\ref{NewtonIterations}) make sense for numerous other choices of the function~$\varphi(\cdot)$,
for instance, for $\varphi(\zeta)=\log \zeta,$ for polynomial and rational~$\varphi(\zeta),$~etc.

The next statement is well-known and is included only for the sake of completeness.

\begin{lemma}
\label{lemma:limitOfNewtonIterationsIsTheInverse}
Suppose that the sequence of Newton's iterations~(\ref{NewtonIterations}) is pointwise convergent in a domain $D\subset\C^n.$
Then its limit $\nu(f):=\lim\limits_{j\rightarrow\infty}\nu_{j}[A,\varphi](f)$ inverts the mapping $f : x \mapsto x + \varphi(Ax)$
as long as the compositions $f\circ\nu$ and $\nu\circ f$ are well-defined.
\end{lemma}
\begin{proof}
We~employ the standard formal argument to show that $f(\nu(f))=f$ and $\nu(f(x))=x$ whenever the compositions are well-defined.
Since both mapping~$f$ and function~$\varphi(\cdot)$ are continuous in their respective domains of definition, it follows that
$$
f(\nu(f))=
f\left(\lim_{j\rightarrow\infty} \nu_{j}(f)\right) =
\lim_{j\rightarrow\infty} \nu_{j}(f) + \varphi\left( A\lim_{j\rightarrow\infty}\nu_{j}(f) \right) =
\lim_{j\rightarrow\infty} \nu_{j}(f) + \lim_{j\rightarrow\infty} \varphi(A\nu_{j}(f))
$$
$$
=\lim_{j\rightarrow\infty} \nu_{j}(f) + \lim_{j\rightarrow\infty}(f - f + \varphi(A\nu_{j}(f))) =
\lim_{j\rightarrow\infty} \nu_{j}(f) + f - \lim_{j\rightarrow\infty} \nu_{j+1}(f) = f.
$$
Similarly $\nu(f(x))=x$ and the lemma follows.
\end{proof}

\begin{definition}\rm
\label{def:finiteNewtonInverse}
If~the sequence of mappings~(\ref{NewtonIterations}) stabilizes, i.e., if there exists $k\in\N$
such that $\nu_{j}=\nu_k$ for all $j\geq k$ then we call its limit~$\nu_k$
{\it the finite Newton's inverse of the mapping $f = x + \varphi(A x)$} and denote it by $\nu_{\rm fin}(f)=\nu_{\rm fin} [A,\varphi](f).$
The smallest integer~$k$ with this property will be called {\it the Newton order of the inverse mapping to~$f$}
and denoted by $\mathfrak{N}(f)=\mathfrak{N}(A,\varphi).$
If~for given~$A$ and~$\varphi$ there is no such finite~$k,$ we will adopt the convention that $\mathfrak{N}(A,\varphi)=\infty.$
\end{definition}

\begin{theorem}
\label{thm:NewtonInverseOfJMap}
For any universal matrix~$U$ and a generic analytic function~$\varphi(\cdot)$
the inverse of the mapping $f=x+\varphi(Ux)$ is given by its finite Newton's inverse $\nu_{\rm fin}[U,\varphi](f)$ and
$\mathfrak{N}(U,\varphi)={\rm rank}\,U.$
\end{theorem}
\begin{proof}
We~denote by $U_1,\ldots,U_n$ the rows of~$U$ and define~$m$ to be the number of different rows in~$U.$
By~Theorem~\ref{thm:UnivMatrixIsDefinedByPartitionAndPermutation} there exist an integer partition $p=(p_1,\ldots,p_m)$
of the dimension~$n$ and a permutation matrix~$\Pi$ of size~$m$ such that the matrix~$U$
is permutation-similar to a matrix of the form~(\ref{canonicalFormOfUMatrix}).
It~suffices to consider the case when the matrix~$U$ is given in its ordered form and when the matrix~$\Pi$ is the identity matrix
since all other cases can be reduced to this one by means of the action of suitable permutation similarity transformations on the matrices~$U$ and~$\Pi.$

We~use induction on the number~$m$ of different rows in the matrix~$U.$
If~$m=1$ then by~Theorem~\ref{thm:UnivMatrixIsDefinedByPartitionAndPermutation} there exist $a_{1}, a_{2}, \ldots, a_{n-1}\in\C$
such that the matrix~$U$ is the $n\times n$ matrix all of whose rows are equal to the vector
$(a_{1}, a_{2}, \ldots, a_{n-1}, -a_{1}-a_{2}-\ldots - a_{n-1}).$
The equations
$$
f_j = x_j + \varphi(a_{1}x_{1} + a_{2}x_{2} + \ldots + a_{n-1}x_{n-1} - (a_{1}+a_{2}+\ldots + a_{n-1})x_{n}), \quad j=1,\ldots,n
$$
imply that
$$
a_{1}f_{1} + a_{2}f_{2} + \ldots + a_{n-1}f_{n-1} - (a_{1}+a_{2}+\ldots + a_{n-1})f_{n} =
$$
$$
a_{1}x_{1} + a_{2}x_{2} + \ldots + a_{n-1}x_{n-1} - (a_{1}+a_{2}+\ldots + a_{n-1})x_{n}
$$
and hence the inverse mapping is given by the first element in the sequence of Newton's iterations~(\ref{NewtonIterations}), i.e.,
$$
x_j = f_j - \varphi(a_{1}f_{1} + a_{2}f_{2} + \ldots + a_{n-1}f_{n-1} - (a_{1}+a_{2}+\ldots + a_{n-1})f_{n}), \quad j=1,\ldots,n,
$$
which concludes the proof in the case when $m = {\rm rank}\,U = 1.$

Let now $m>1$ and denote $\tilde{p}:=(p_1,\ldots,p_{m-1}),$ $|\tilde{p}|:=p_1+\ldots+p_{m-1},$ and $\tilde{x}=(x_1,\ldots,x_{|\tilde{p}|}).$
Since the universal matrix~$U$ is assumed to be given in its ordered form and the matrix~$\Pi$ is the identity matrix,
it follows by Theorem~\ref{thm:UnivMatrixIsDefinedByPartitionAndPermutation} that the submatrix~$\tilde{U}$ of the matrix~$U$
comprising the first~$|\tilde{p}|$ rows and first~$|\tilde{p}|$ columns of~$U$ is universal, too.
By~the induction hypothesis the conclusion of the theorem holds for the mapping $\tilde{f}:=\tilde{x}+\varphi(\tilde{U}\tilde{x}).$

Let now~$\hat{U}$ denote the $|\tilde{p}|\times |p|$-submatrix of the matrix~$U$ comprising the first~$|\tilde{p}|$ rows of~$U.$
Composing the mapping~$\tilde{f}$ with the shift $\tilde{x}\mapsto\tilde{x}+c$ for a suitable $c\in\C^{|\tilde{p}|}$ we conclude that
the mapping $\hat{f}:\tilde{x}\mapsto \tilde{x} + \varphi(\hat{U} x)$ viewed as a mapping with the variables $x_{1},\ldots,x_{|\tilde{p}|}$
depending on the parameters $x_{|\tilde{p}|+1},\ldots,x_{|p|}$ has unit Jacobian $J(\hat{f}; \tilde{x}).$
Furthermore the sequence of Newton's iterations~(\ref{NewtonIterations}) yields the inverse mapping to~$\hat{f}$
and $\mathfrak{N}(\tilde{f})=\mathfrak{N}(\hat{f}).$

Since by assumption the universal matrix $U=(u_{jk})$ is given in its ordered form and the matrix~$\Pi$ is the identity matrix,
the last~$p_m$ columns of the matrix~$U$ comprise block submatrices of type~$Z$
(i.e.,~each of these blocks consists of equal rows whose elements sum up to zero).
Thus $u_{|p|,|\tilde{p}|+1} + \ldots + u_{|p|,|p|} = 0$ and the last~$p_m$ equations in the system $f=x+\varphi(Ux)$ imply that
$$
u_{|p|,|\tilde{p}|+1} x_{|\tilde{p}|+1} + \ldots + u_{|p|,|p|} x_{|p|}  = u_{|p|,|\tilde{p}|+1} f_{|\tilde{p}|+1} + \ldots + u_{|p|,|p|} f_{|p|}.
$$
Furthermore $f_j - f_k = x_j - x_k$ for any $j,k=|\tilde{p}|+1,\ldots,|p|$ and hence any linear combination of $f_j,\,j=|\tilde{p}|+1,\ldots,|p|$
whose coefficients sum up to zero is also a linear combination of $x_j,\,j=|\tilde{p}|+1,\ldots,|p|.$
It~follows that one more iteration~(\ref{NewtonIterations}) applied to the solution of the system $\hat{f}=\tilde{x}+\varphi(\hat{U}x)$
yields the inverse to $x+\varphi(Ux),$ which concludes the proof.
\end{proof}
\begin{example}
\label{ex:universalMatrixInDim=5}\rm
For any $a,b,c,s,t,u,v\in\C$ the matrix
$$
M=
\left(
\begin{array}{ccccc}
 a & -a & b & c & -b-c \\
 a & -a & b & c & -b-c \\
 s &  t & u & v & -u-v \\
 s &  t & u & v & -u-v \\
 s &  t & u & v & -u-v \\
\end{array}
\right)
$$
is universal with generic rank~2. It~is defined by the partition $2+3$ of the dimension $n=5$ together with the trivial permutation on two elements.
For $x\in\C^5$ and an analytic function~$\varphi(\cdot)$ the inverse $x=x(f)$
to the mapping $f:x \mapsto x+\varphi(M x)$ is given by the second iteration of the Newton mapping~(\ref{NewtonIterations}):
$$
{\small
\begin{array}{c}
f_1-\varphi\left(a (f_1-f_2) + b(f_3-f_5) + c(f_4-f_5) \right), \\
f_2-\varphi\left(a (f_1-f_2) + b(f_3-f_5) + c(f_4-f_5) \right),\\
f_3-\varphi\left(s f_1 + t f_2 + u(f_3-f_5) + v(f_4-f_5) - (s+t) \varphi\left(a (f_1-f_2) + b(f_3-f_5) + c(f_4-f_5)\right)\right), \\
f_4-\varphi\left(s f_1 + t f_2 + u(f_3-f_5) + v(f_4-f_5) - (s+t) \varphi\left(a (f_1-f_2) + b(f_3-f_5) + c(f_4-f_5)\right)\right), \\
f_5-\varphi\left(s f_1 + t f_2 + u(f_3-f_5) + v(f_4-f_5) - (s+t) \varphi\left(a (f_1-f_2) + b(f_3-f_5) + c(f_4-f_5)\right)\right). \\
\end{array}
}
$$
Straightforward computation shows that $f\circ x \equiv x\circ f \equiv {\rm Id}$ as long as both compositions are well-defined.
\end{example}

Theorem~\ref{thm:NewtonInverseOfJMap} together with the explicit form of the Newton mapping~(\ref{NewtonIterations}) yield the following corollary.

\begin{corollary}
\label{cor:inverseOfJMapDefinedByUMatrixIsComposition}
Let~$U$ be a universal matrix with the rows $U_1,\ldots,U_n$ and let~$\varphi(\cdot)$ be a univariate analytic function
defined in a nonempty domain in the complex plane.
The components of the inverse to the mapping $x+\varphi(Ux)$ (that is, the solution to the system of equations
$f_j=x_j+\varphi(\langle U_j,x \rangle),$ $j=1,\ldots,n$) are finite superpositions of multiplication by a constant,
addition, and the function~$\varphi(\cdot)$ applied to $f_1,\ldots,f_n.$
\end{corollary}

Due to the classical formulation of the Jacobian Conjecture the special case when the function~$\varphi(\cdot)$
is a polynomial is of particular importance.

\begin{corollary}
\label{cor:inverseOfJMapDefinedByUMatrixAndPolynomialPhiIsPolynomial}
The inverse of the mapping $x+\varphi(Ux)$ defined by a universal matrix~$U$ and a polynomial~$\varphi$ is a polynomial mapping, too.
\end{corollary}


\section{Jacobian equations for mappings of the form $x + (Ax)^d$
\label{sec:JEquationsFor x+(Ax)^d}}

In this section we restrict our attention to the crucially important special case when the mapping $x+\varphi(Ax)$
is defined by a monomial function $\varphi(\zeta)=\zeta^d.$
That is, we consider polynomial mappings of the form $x + (Ax)^d,$ where $x=(x_1,\ldots,$ $x_n)$ $\in\C^n,$
$A=(a_{jk})$~is $n\times n$-matrix, and~$d\in\N.$
The components of this mapping are given~by
\begin{equation}
\label{JUniformization}
x_j + \left(\sum_{k=1}^n a_{jk}x_k\right)^d, \quad j=1,\ldots,n.
\end{equation}
If~$d=0$ or $d=1$ then the mapping~(\ref{JUniformization}) is (affine)
linear and its Jacobian is the determinant of the corresponding matrix.
We~do not consider these trivial cases.

It~follows from~(\ref{universalMatrixEquation}) that the mapping~(\ref{JUniformization}) has unit Jacobian if and only if the elements of the matrix~$A$
satisfy a system of algebraic equations which depends on the integers~$n$ and~$d.$
We~adopt the following definition.

\begin{definition}\rm
\label{def:JEquations}
By~the {\it Jacobian equations in dimension~$n\geq 2$ of degree~$d,$ $d\in\N$} we will mean the algebraic equations
with the variables~$a_{jk}$ which define the matrices~$A=(a_{jk})$ such that the Jacobian of the mapping
$x + (Ax)^d$ is identically equal to~$1.$
\end{definition}

\begin{example}\rm
\label{ex:JEquationsDim=2,Deg=2}
The simplest nontrivial case is that of a bivariate mapping of degree two, i.e., the mapping of the form
$$
\begin{array}{c}
f_1 = x_1 + (a_{11} x_1 + a_{12} x_2)^2, \\
f_2 = x_2 + (a_{21} x_1 + a_{22} x_2)^2.
\end{array}
$$
The Jacobian of this mapping is the following polynomial of degree two in $x_1,x_2$ whose coefficients involve the determinant
and the permanent of the matrix $A=(a_{jk}),$ $j,k=1,2:$
$$
\begin{array}{r}
J(f_1,f_2; x_1,x_2) =
1
+2 x_1 \left(a_{11}^2+a_{21} a_{22}\right)
+2 x_2 \left(a_{11} a_{12}+a_{22}^2\right)
+4 x_1^2 a_{11} a_{21} \left(a_{11} a_{22}-a_{12} a_{21}\right)\phantom{.} \\
+4 x_1 x_2 \left(a_{11} a_{22}-a_{12} a_{21}\right) \left(a_{11} a_{22}+a_{12} a_{21}\right)
+4 x_2^2 a_{12} a_{22} \left(a_{11} a_{22}-a_{12} a_{21}\right).
\end{array}
$$
Thus the Jacobian equations in dimension~2 and of degree~2 constitute the following system of algebraic equations:
\begin{equation}
\label{JEquations,dim=2,deg=2}
\left\{
\begin{array}{l}
a_{11}^2 + a_{21} a_{22} = 0, \\
a_{11} a_{12} + a_{22}^2 = 0, \\
a_{12} a_{22}(a_{11} a_{22} - a_{12} a_{21}) = 0, \\
a_{11} a_{21}(a_{11} a_{22} - a_{12} a_{21}) = 0, \\
(a_{11} a_{22} - a_{12} a_{21}) (a_{11} a_{22} + a_{12} a_{21}) = 0.\\
\end{array}
\right.
\end{equation}
This system can be simplified and the number of equations can be reduced.
\end{example}

The number and complexity of Jacobian equations rapidly grow with the dimension of the ambient space and the degree of the mapping.

\begin{example}\rm
\label{ex:dim=3,deg=2,JEquations}
The Jacobian equations for the degree two mapping in three-dimensional space are given by
$$
\left\{
\begin{array}{l}
a_{11}^2 + a_{21} a_{22} + a_{31} a_{33} = 0, \\
a_{11} a_{12} + a_{22}^2 + a_{32} a_{33} = 0, \\
a_{11} a_{13} + a_{22} a_{23} + a_{33}^2 = 0, \\
a_{21} a_{22} a_{11}^2+a_{31} a_{33} a_{11}^2-a_{13} a_{31}^2 a_{11}-a_{11} a_{12} a_{21}^2+a_{21} a_{22} a_{31} a_{33}-a_{21} a_{23} a_{31} a_{32} = 0, \\
a_{21} a_{22} a_{12}^2+a_{22} a_{23} a_{32}^2-a_{11} a_{22}^2 a_{12}-a_{22}^2 a_{32} a_{33}+a_{13} a_{31} a_{32} a_{12}-a_{11} a_{32} a_{33} a_{12} = 0, \\
a_{31} a_{33} a_{13}^2+a_{23}^2 a_{32} a_{33}-a_{22} a_{23} a_{33}^2-a_{11} a_{33}^2 a_{13}+a_{12} a_{21} a_{23} a_{13}-a_{11} a_{22} a_{23} a_{13} = 0, \\
a_{11}^2a_{22}^2 - a_{12}^2a_{21}^2 + a_{11}^2a_{32}a_{33} + a_{22}^2a_{31}a_{33} - a_{12}a_{13}a_{31}^2 - a_{21}a_{23}a_{32}^2 + \\
a_{11}a_{12}a_{31}a_{33} + a_{21}a_{22}a_{32}a_{33} - a_{11}a_{13}a_{31}a_{32} - a_{22}a_{23}a_{31}a_{32}  = 0, \\
a_{11}^2a_{33}^2 - a_{13}^2a_{31}^2 + a_{11}^2a_{22}a_{23} + a_{21}a_{22}a_{33}^2 - a_{23}^2a_{31}a_{32}  - a_{12}a_{13}a_{21}^2 + \\
a_{11}a_{13}a_{21}a_{22} + a_{22}a_{23}a_{31}a_{33} - a_{11}a_{12}a_{21}a_{23} - a_{21}a_{23}a_{32}a_{33} = 0, \\
a_{22}^2a_{33}^2 - a_{23}^2a_{32}^2 + a_{11}a_{13}a_{22}^2 + a_{11}a_{33}^2a_{12} - a_{13}^2a_{31}a_{32} - a_{21}a_{23}a_{12}^2 + \\
a_{11}a_{22}a_{23}a_{12} + a_{11}a_{13}a_{32}a_{33} - a_{13}a_{21}a_{22}a_{12}  - a_{13}a_{31}a_{33}a_{12} = 0, \\
a_{11} a_{21} a_{31} \det A = a_{12} a_{22} a_{32} \det A = a_{13} a_{23} a_{33} \det A = 0, \\
(a_{11} a_{22} a_{31} + a_{11} a_{21} a_{32} + a_{12} a_{21} a_{31}) \det A = 0, \\
(a_{11} a_{23} a_{31} + a_{11} a_{21} a_{33} + a_{13} a_{21} a_{31}) \det A = 0, \\
(a_{11} a_{22} a_{32} + a_{12} a_{22} a_{31} + a_{12} a_{21} a_{32}) \det A = 0, \\
(a_{11} a_{23} a_{33} + a_{13} a_{23} a_{31} + a_{13} a_{21} a_{33}) \det A = 0, \\
(a_{12} a_{23} a_{32} + a_{12} a_{22} a_{33} + a_{13} a_{22} a_{32}) \det A = 0, \\
(a_{12} a_{23} a_{33} + a_{13} a_{23} a_{32} + a_{13} a_{22} a_{33}) \det A = 0, \\
\det A \, \rm{perm}\, A = 0. \\
\end{array}
\right.
$$
The system of Jacobian equations defined by cubic mapping in dimension four comprises~294 equations.
The highest of the total degrees of these equations with respect to the variables~$a_{jk}$ equals~48.
\end{example}

Despite the formidable complexity of the general system of Jacobian equations, it contains a subsystem whose structure is fairly transparent.
This subsystem comprises the first two equations in~(\ref{JEquations,dim=2,deg=2}).
For the Jacobian equations in Example~\ref{ex:dim=3,deg=2,JEquations}, this subsystem consists of the first three equations.
We~adopt the next definition.

\begin{definition}
\label{def:HadamardProduct}\rm
Let $A=(a_{jk})$ and $B=(b_{jk})$ be matrices of equal size. By~$A\odot B$ we will denote the Hadamard (termwise) product
of~$A$ and~$B,$ i.e., the matrix~$(a_{jk}b_{jk}).$ For a square matrix~$A$ and~$d\in\N$ we will denote by~$A^{\odot d}$
the $d$th Hadamard power of~$A,$ i.e., the matrix $A\odot\ldots\odot A$ ($d$~copies of~$A$).
\end{definition}

The following lemma holds.

\begin{lemma}
\label{lemma:simpleJEquations}
Let $A=(a_{jk})$ be a square matrix of size $n\geq 2.$ For any $d=2,3,\ldots$ the system of Jacobian equations in dimension~$n$ of degree~$d$
contains the {\rm simple Jacobian equations}
\begin{equation}
\label{simpleJEquations}
\left(A^{T}\right)^{\odot(d-1)} \, {\rm diag}\, A = 0.
\end{equation}
Here~$A^{T}$ denotes the transpose of~$A$ and ${\rm diag\, A}$ is the vector of diagonal elements of~$A.$
\end{lemma}
\begin{proof}
Induction on~$n$ shows that the components of the left-hand side of~(\ref{simpleJEquations})
are the coefficients by the monomials~$x_{j}^{d-1},$ $j=1,\ldots,n$ in the Jacobian of the mapping $x+(Ax)^d.$
Since these monomials are linearly independent, the mapping can only be Jacobian if all of the corresponding coefficients vanish simultaneously.
\end{proof}

The above examples suggest that the Jacobian equations are homogeneous with respect to the variables~$a_{jk}.$
This is indeed the case and it is moreover possible to explicitly describe all homogeneities of the Jacobian equations in any dimension and degree.
The following lemma holds.

\begin{lemma}
\label{lemma:homogeneitiesOfJEquations}
For $s=(s_1,\ldots,s_n)\in\left(\C^*\right)^n$ define $H_{n,d}(s)$ to be the $n\times n$-matrix whose $j$th row equals
$s_1\ldots [j] \ldots s_n \cdot \left(s_1^d, s_2^d,\ldots,s_n^d\right),$ $j=1,\ldots,n.$
For $v=(v_1,\ldots,v_n)\in\C^n$ we adopt the notation $s^v := s_{1}^{v_1} \ldots s_{n}^{v_n}.$
This matrix $H_{n,d}(s)$ is explicitly given by
\begin{equation}
\label{homogeneitiesOfJEquations}
H_{n,d}(s) = (s^{\mathcal{I} - e_j + d e_k}), \quad j,k = 1,\ldots,n,
\end{equation}
where $\mathcal{I}=(1,\ldots,1)$ is the vector of~$n$ units, $\{ e_j = (0,\ldots,1,\ldots,0) \}$
is the standard basis in the $n$-dimensional space, $\delta_{jk}$ is the Kronecker delta.
The Jacobian equations in dimension~$n$ for mapping of degree~$d$ are homogeneous with respect to the matrix~$H_{n,d}(s).$
That is, for any matrix~$A$ whose elements satisfy the Jacobian equations the elements of the Hadamard (termwise) product $A\odot H_{n,d}(s)$
also satisfy these equations.

Moreover, the matrix $H_{n,d}(s)$ represents all homogeneities of the Jacobian equations of degree~$d$ in dimension~$n.$
\end{lemma}
\begin{proof}
Let~$A$ be any square matrix of size~$n.$ It follows from~(\ref{universalMatrixEquation}) that the identity
$$
J\left(x+((H_{n,d}(s)\odot A) x)^d; x\right) \equiv J\left(\xi + (A\xi)^d; \xi\right)
\bigg|_{\xi_j = x_j s_{j}^d\left(\prod\limits_{k=1}^n s_k\right)^\frac{d}{d-1}, \,\, j=1,\ldots,n}
$$
holds for any $d=2,3,\ldots$ and any $x\in\C^n.$ Since the Jacobian equations are defined through the coefficients by the linearly independent
monomials in $J(\xi + (A\xi)^d; \xi),$ it follows that the matrix~$A$ forms a good pair with the function~$\zeta^d$ if and only if
so does the Hadamard product $H_{n,d}(s)\odot A.$

The Jacobian equations of degree~$d$ in dimension~$n$ do not admit any homogeneities other than those encoded by the matrix~$H_{n,d}(s)$
since this is already true for the simple Jacobian equations~(\ref{simpleJEquations}).
\end{proof}

Despite the complex structure of the general system of Jacobian equations of a given degree~$d$,
in the bivariate case it turns out to be possible to explicitly solve it by parameterizing its set of solutions.

\begin{lemma}
\label{lemma:parametrizationOfJMatrices,dim=2}
Let $d\in\C,$ $d\neq 0,1.$
A~nonzero $2\times 2$-matrix forms a good pair with the function~$\zeta^d$ if and only if
there exist $(s,t)\in\P^1$ such that this matrix has the form
\begin{equation}
\label{parametrizationOfJMatrices,dim=2}
G(2,d):=
\left(
\begin{array}{cc}
s t^d   & -s^{d+1} \\
t^{d+1} & -s^d t   \\
\end{array}
\right)
\equiv
\left(
\begin{array}{cc}
1 & -1 \\
1 & -1 \\
\end{array}
\right)\odot
\left(
\begin{array}{cc}
s t^d   & s^{d+1} \\
t^{d+1} & s^d t   \\
\end{array}
\right).
\end{equation}

The set of solutions to the Jacobian equations of degree~$d$ in the bivariate case has dimension~$1$ regardless of the value of $d\neq 0,1.$
Furthermore, the inverse of the mapping $f = x+(G(2,d)x)^d$ is given by the first iteration of the Newton mapping~(\ref{NewtonIterations}),
i.e., by $x = f - (G(2,d)f)^d.$
\end{lemma}

We~remark that $G(2,d)$ is the Hadamard product of the base of the double cone of nontrivial (i.e.,~with no zero rows or columns)
universal matrices of size~$2$ and the matrix of homogeneities~$H_{2,d}(s,t).$
Besides, it is easy to see that all the three universal matrices in Example~\ref{ex:unversalMatricesInDim2} are special instances
of the matrix $G(2,d)$ corresponding to $(s,t)=(\sqrt[d+1]{-a},0),$ $(s,t)=(0,\sqrt[d+1]{a}),$ and $(s,t)=(\sqrt[d+1]{a},\sqrt[d+1]{a}),$
respectively.

\begin{proof}
Straightforward calculation shows that the matrix~(\ref{parametrizationOfJMatrices,dim=2}) indeed forms a good pair with the function~$\zeta^{d},$
i.e., that the Jacobian of the mapping $f_{1}=x_{1} + (s t^d x_1 - s^{d+1} x_2)^d,$ $f_{2}=x_{2} + (t^{d+1} x_1 - s^d t x_2)^d$ is identically equal to~1.

Conversely, let the matrix
$
A=
\left(
\begin{array}{cc}
a_{11} & a_{12} \\
a_{21} & a_{22} \\
\end{array}
\right)
$
form a good pair with the function~$\zeta^{d},$ $d\neq 0,1.$
By~(\ref{universalMatrixEquation}) the Jacobian of the mapping $f_{j} = x_{j} + (a_{j1}x_1 + a_{j2}x_2)^d,$ $j=1,2$ is given by
\begin{equation}
\label{JInDim=2,Deg=d}
\begin{array}{cc}
1 + d a_{11}( a_{11} x_1 + a_{12}x_2 )^{d-1} + d a_{22}( a_{21} x_1 + a_{22}x_2 )^{d-1} + \\
d^2 (a_{11}a_{22} - a_{12}a_{21}) ( a_{11} x_1 + a_{12}x_2 )^{2d-2} ( a_{21} x_1 + a_{22}x_2 )^{2d-2}.
\end{array}
\end{equation}
Since $d\neq 1,$ the homogeneous components of~(\ref{JInDim=2,Deg=d}) of degrees $d-1$ and $2d-2$ in the variables~$x_1,x_2$ must both vanish
for the mapping $x+(Ax)^d$ to be Jacobian. In~particular, $\det A = 0.$
If~$A$ is the zero matrix then the conclusion of the lemma is trivially true.
Thus we may without loss of generality assume that the second row of~$A$ contains a nonzero element, say,~$a_{21}.$
Since~$A$ is degenerate, it follows that $(a_{11},a_{12})=\lambda(a_{21},a_{22})$ for some $\lambda\in\C.$
We~therefore have $a_{11}=\lambda a_{21},$ $a_{22}=-\lambda^d a_{21},$ and $a_{12}=-\lambda^{d+1}a_{21}.$
Setting without loss of generality $a_{21}=1$ we conclude that
$$
A=
\left(
\begin{array}{cl}
\lambda & -\lambda^{d+1} \\
  1     & -\lambda^d     \\
\end{array}
\right).
$$
Passing to the projective coordinates $(s,t)\in\P^{1}$ such that $\lambda=s/t$ we arrive at~(\ref{parametrizationOfJMatrices,dim=2}).

It~remains to check the inverse mapping. The composition $x_{1}(f_{1}(x_1,x_2), f_{2}(x_1,x_2))$ is given by
$$
\begin{array}{c}
f_1 - (s t^d f_{1} - s^{d+1} f_2)^d = \\
x_1 + \left(s t^d x_{1} - s^{d+1} x_2\right)^d - \\
\left(s t^d (x_1 + (s t^d x_{1} - s^{d+1} x_2)^d) - s^{d+1} (x_{2} + (t^{d+1} x_1 - s^d t x_2)^d)\right)^d \equiv x_1 \\
\end{array}
$$
by power series expansion. Similarly the composition $x_{2}(f_{1}(x_1,x_2), f_{2}(x_1,x_2))$ is identically equal to~$x_2$
and furthermore $f(x_{1}(f_1,f_2), x_{2}(f_1,f_2))\equiv f.$ The proof is complete.
\end{proof}

Using Theorem~\ref{thm:UnivMatrixIsDefinedByPartitionAndPermutation}, Lemma~\ref{homogeneitiesOfJEquations},
and Corollary~\ref{cor:inverseOfJMapDefinedByUMatrixAndPolynomialPhiIsPolynomial},
we arrive at the next theorem which gives a family of matrices that form good pairs with the monomial function of prescribed degree~$d.$

\begin{theorem}
\label{thm:HadamardProductsOfUniversalAndHomogeneitiesAreGood}
Let $n\geq 2$.
For any partition $p=(p_1,\ldots,p_m)$ of~$n,$ any permutation matrix of size~$m,$ and any $s\in\C^n$
the mapping $f:\C^n\rightarrow\C^n$ defined through
$$
f : x \mapsto x + \left((U(p,\Pi)\odot H_{n,d}(s))\,x\right)^d
$$
has unit Jacobian and is polynomially invertible. The inverse mapping is given by the finite Newton's inverse $\nu_{\rm fin}(f)$ and $\mathfrak{N}(f)=m.$
\end{theorem}

One could hope that the space of all matrices that form good pairs with at least one monomial function $\zeta^d,\, d=2,3,4,\ldots$
admits a base comprising all universal matrices and is fibrated by taking the Hadamard product with the matrix of homogeneities $H_{n,d}(s)$
of the Jacobian equations. Lemma~\ref{lemma:parametrizationOfJMatrices,dim=2} yields that this is indeed the case in two dimensions.
The next example shows that this statement, unfortunately, does not hold in general.

\begin{example}
\label{ex:goodMatr,dim=4,deg=2,rank=2,centrosymmetric,idempotent}\rm
For any $a,b,c\in\C$ such that $ab(a+b)(2ac+b^2)\neq 0$ the idempotent centrosymmetric matrix
$$
{\tiny
M=
\frac{1}{a b (a+b) \left(2 a c+b^2\right)}
\left(
\begin{array}{cccc}
 a^2 b \left(b^2+2 a c\right) & -b \left(b^2+a c\right) \left(b^2+2 a c\right) & a b c \left(b^2+2 a c\right) & 0 \\
 -a^3 \left(b^2+a c\right) & a b^2 \left(b^2+2 a c\right) & 0 & -a^4 c \\
 -a^4 c & 0 & a b^2 \left(b^2+2 a c\right) & -a^3 \left(b^2+a c\right) \\
 0 & a b c \left(b^2+2 a c\right) & -b \left(b^2+a c\right) \left(b^2+2 a c\right) & a^2 b \left(b^2+2 a c\right) \\
\end{array}
\right)
}
$$
forms a good pair with the function~$\zeta^2.$
The matrix~$M$ is {\it not} of the form $U\odot H$ for any universal matrix~$U$
and any choice of $s\in\C^4$ in the matrix of homogeneities $H=H_{4,2}(s).$
Indeed, any such product either consists of nonzero elements only,
or contains a certain positive number of zero elements in its diagonal, which is not the case for the matrix~$M.$
The solution to the system of equations $f=x+(Mx)^2$ is given by
the second iteration of the Newton mapping~(\ref{NewtonIterations}), i.e., $x=f-(M(f-(Mf)^2))^2.$
\end{example}
\begin{theorem}
\label{thm:HadamardProductParametrizationForGoodMatrixWithProportionalRows}
Let~$M$ be a square matrix of size~$n$ such that

1)~The mapping $f=x+(Mx)^d$ has unit Jacobian;

2)~Any set of linearly dependent rows of~$M$ contains proportional rows.
In~other words, any linear relation between the rows of~$M$ involves at most two rows which are thereby proportional.

Then there exist a universal matrix~$U$ and $s\in\C^n$ such that $M=U\odot H_{n,d}(s).$
Furthermore, the corresponding Jacobian mapping $x+(Mx)^d$ is polynomially invertible, its inverse is given by $\nu_{\rm fin}(f),$
and $\mathfrak{N}(f)={\rm rank}\,M.$
\end{theorem}
\begin{proof}
We~employ the same idea as the one used in the proof of Theorem~\ref{thm:NewtonInverseOfJMap}.
The number of different rows of a universal matrix is the counterpart of the number of
linearly independent rows of a matrix that forms a good pair with a monomial function,
so we argue by induction on ${\rm rank}\,M.$

To~establish the base of induction, we need to provide a decomposition of a general matrix of rank~1
that forms a good pair with $\varphi(\zeta)=\zeta^d$ into the Hadamard product of a universal matrix and a matrix of homogeneities.
The rows of a rank~1 matrix~$M$ of size~$n$ are all proportional to a certain single vector, say, $v=(v_1,\ldots,v_n)\in\C^n.$
Thus a general rank~1 matrix can be written in the form
\begin{equation}
\label{generalRankOneMatrix}
M=
\left(
\begin{array}{cccc}
c_1 v_1 & c_1 v_2 & \ldots & c_1 v_n \\
c_2 v_1 & c_2 v_2 & \ldots & c_2 v_n \\
 \ldots & \ldots  & \ldots & \ldots  \\
c_n v_1 & c_n v_2 & \ldots & c_n v_n \\
\end{array}
\right) = c\otimes v,
\end{equation}
where $c=(c_1,\ldots,c_n)\in\C^n$ is the vector of the coefficients of proportionality.

It~follows immediately from~(\ref{universalMatrixEquation}) that the matrix~(\ref{generalRankOneMatrix})
forms a good pair with the function~$\zeta^d,$ $d\in\C^{*}$ if and only if
\begin{equation}
\label{whenRankOneMatrixFormsAGoodPairWithPowerd}
c_{1}^{d} v_1 + c_{2}^{d} v_2 + \ldots + c_{n}^{d} v_n = 0.
\end{equation}
Solving this equation for $v_n = -\sum\limits_{j=1}^{n-1} \left(\frac{c_j}{c_n}\right)^d v_j,$
introducing new parameters $a_{j}^{(1,1)} := \left(\frac{c_j}{c_n}\right)^d v_j,$ $j=1,\ldots,n-1,$ and
$s_j:=1/c_j,$ $j=1,\ldots,n,$ and clearing denominators,
we conclude that the matrix~$M$ forms a good pair with the function~$\zeta^d$ if and only if
$$
M=
\left(
\begin{array}{ccccc}
s_1^d   s_2 s_3 s_4 \ldots s_n  \, a^{(1,1)}_1 & s_2^{d+1} s_3 s_4 \ldots s_n     \, a^{(1,1)}_2 & \ldots  & -s_2 s_3 s_4 \ldots s_n^{d+1} \sum\limits_{j=1}^{n-1} a^{(1,1)}_j \\
s_1^{d+1} s_3 s_4 \ldots s_n    \, a^{(1,1)}_1 & s_1 s_2^d s_3 s_4 \ldots s_n     \, a^{(1,1)}_2 & \ldots  & -s_1 s_3 s_4 \ldots s_n^{d+1} \sum\limits_{j=1}^{n-1} a^{(1,1)}_j \\
            \ldots                             & \ldots                                          & \ldots  & \ldots                                              \\
s_1^{d+1} s_2 s_3 \ldots s_n    \, a^{(1,1)}_1 & s_1 s_2^{d+1} s_3 \ldots s_n     \, a^{(1,1)}_2 & \ldots  & -s_1 s_2 s_3 \ldots s_n^{d+1} \sum\limits_{j=1}^{n-1} a^{(1,1)}_j \\
s_1^{d+1} s_2 s_3 \ldots s_{n-1}\, a^{(1,1)}_1 & s_1 s_2^{d+1} s_3 \ldots s_{n-1} \, a^{(1,1)}_2 & \ldots  & -s_1 s_2 s_3 \ldots s_n^d \sum\limits_{j=1}^{n-1} a^{(1,1)}_j \\
\end{array}
\right)=
$$
$$
= U(n,\varepsilon_1) \odot H_{n,d}(s).
$$
That is, any matrix of rank~1 that forms a good pair with the function~$\zeta^d$ can be represented as the Hadamard product
of the matrix~(\ref{homogeneitiesOfJEquations}) of homogeneities of the Jacobian equations in dimension~$n$ and degree~$d$ and
the universal matrix $U(n,\varepsilon_1)$ (see~(\ref{canonicalFormOfUMatrix})) defined by the trivial partition of the dimension~$n$
and the trivial permutation~$\varepsilon_1$ of its single element.

If~the matrix~$M$ satisfies the condition~(\ref{whenRankOneMatrixFormsAGoodPairWithPowerd}) then the inverse mapping,
i.e., the solution to the system of equations
$f_{j}=x_{j} + (c_j v_1 x_1 + c_j v_2 x_2 + \ldots + c_j v_n x_n)^d,$ $j=1,\ldots,n$
is given~by
$$
x_{j} = f_{j} - (c_j v_1 f_1 + c_j v_2 f_2 + \ldots + c_j v_n f_n)^d, \quad j=1,\ldots,n.
$$
Using the convention that the function~$\zeta^d$ acts termwise on the coordinates of its vector argument,
we can write the above inverse in the concise form $x = f - (Mf)^d.$

The inductive step is completely parallel to the one used in the proof of Theorem~\ref{thm:NewtonInverseOfJMap} and we omit~it.
\end{proof}

The next example shows that the linear relations between the rows of a matrix that forms a good pair with a monomial function
in general need not be proportionalities of the rows and hence Theorem~\ref{thm:HadamardProductParametrizationForGoodMatrixWithProportionalRows}
does not apply to all matrices that form good pairs with monomial functions.

\begin{example}
\label{ex:dim=6,deg=3,rank=2,two-dim family} \rm
For any $a,b\in\C$ the matrix
\begin{equation}
\label{dim=6,deg=3,rank=2,two-dim family}
M(a,b):=
\left(
\begin{array}{cccccc}
 20 a & 0 & 8 a & -4 a & -3 a & a \\
 0 & 20 b & -4 b & 8 b & b & -3 b \\
 40 a & 20 b & 4 (4 a-b) & -8 (a-b) & b-6 a & 2 a-3 b \\
 20 a & 40 b & 8 (a-b) & 4 (4 b-a) & 2 b-3 a & a-6 b \\
 60 a & 20 b & 4 (6 a-b) & 4 (2 b-3 a) & b-9 a & 3 (a-b) \\
 20 a & 60 b & 4 (2 a-3 b) & 4 (6 b-a) & -3 (a-b) & a-9 b \\
\end{array}
\right)
\end{equation}
forms a good pair with the function $\varphi(\zeta)=\zeta^3.$
This matrix has generic rank~$2,$ its last four rows being linear combinations of the first two rows.
The condition~2) in Theorem~\ref{thm:HadamardProductParametrizationForGoodMatrixWithProportionalRows} is violated for the matrix~$M(a,b)$
for generic choice of $a,b\in\C^2.$
A~basis in the two-dimensional linear space $\{ M(a,b): a,b\in\C \}$ is given by the permutation-similar matrices $M(1,0)$ and~$M(0,1).$
The solution to the system of equations $f=x+(M(a,b)x)^3$ is given by the second iteration of the Newton mapping~(\ref{NewtonIterations}).
Similarly to Example~\ref{ex:goodMatr,dim=4,deg=2,rank=2,centrosymmetric,idempotent}, the matrix $M(a,b)$ is not a special instance
of the Hadamard product of a universal matrix and a matrix of homogeneities.

The matrix $M(a,b)$ is a representative of a family of matrices that form good pairs with $\varphi(\zeta)=\zeta^3$ and whose elements depend on eight
additional algebraically independent complex parameters whose values have been set to integer constants in this example.
The generic element of this family is however by far too cumbersome to be displayed in a paper.
\end{example}


\section{Examples}
\label{sec:examples}

Let $A=(a_{jk})$ be a square matrix of size $n\geq 2.$
Consider the mapping $x+\log(Ax)$ whose coordinates are given by
\begin{equation}
\label{logMap}
x_j + \log\left(\sum_{k=1}^n a_{jk}x_k\right), \quad j=1,\ldots,n.
\end{equation}
Due to the main property of the logarithmic function (i.e., $\log(ab)=\log a + \log b$) the set of matrices $A=(a_{jk})$ defining
Jacobian mappings of the form~(\ref{logMap}) admits row-wise homogeneity. That is, if a matrix~$A$ with the rows $A_1,\ldots,A_n$
forms a good pair $(A,\log)$ then for any constants $c_j\in\C^{*},$ $j=1,\ldots,n$
the matrix with the rows $c_1 A_1,\ldots,c_n A_n$ also forms a good pair with the logarithmic function.

\begin{example}
\label{ex:JEquations,LOG,dim=2}\rm
In the bivariate case, the Jacobian equations associated with the mapping~(\ref{logMap}) are as follows:
$$
\left\{
\begin{array}{l}
a_{11} (a_{21}+a_{22})=0,      \\
a_{22} (a_{11}+a_{12})=0,      \\
a_{11} a_{22}-a_{12} a_{21}=0. \\
\end{array}
\right.
$$
Thus any matrix which does not contain zero rows and such that the mapping~(\ref{logMap}) has unit Jacobian is of the form
$$
\left(
\begin{array}{cc}
a & -a \\
b & -b \\
\end{array}
\right),
$$
where $a,b\in\C^{*}.$ Using the main property of the logarithmic function we may without loss of generality choose arbitrary nonzero values of $a,b.$
We~remark that for $a=1,$ $b=-1$ we obtain the circulant matrix $C(1,-1).$
\end{example}

\begin{example}
\label{ex:JEquations,LOG,dim=3}\rm
For $n=3,$ the equations defining the set of matrices
$$
A=
\left(
\begin{array}{cccc}
a_{11} & a_{12} & a_{13} \\
a_{21} & a_{22} & a_{23} \\
a_{31} & a_{32} & a_{33} \\
\end{array}
\right)
$$
such that the Jacobian of the mapping~(\ref{logMap}) is equal to~$1$ are as follows:
$$
\left\{
\begin{array}{l}
a_{11} \left(a_{21} a_{31}+a_{22} a_{31}+a_{21} a_{33}\right) = 0, \\
a_{22} \left(a_{11} a_{32}+a_{12} a_{32}+a_{12} a_{33}\right) = 0, \\
a_{33} \left(a_{13} a_{22}+a_{11} a_{23}+a_{13} a_{23}\right) = 0, \\
a_{11}a_{22}a_{33} + a_{11}a_{22}a_{31} + a_{12}a_{22}a_{31} + a_{11}a_{21}a_{32} + a_{11}a_{22}a_{32} + a_{12}a_{21}a_{33} = 0, \\
a_{13}a_{22}a_{31} + a_{11}a_{23}a_{31} + a_{11}a_{21}a_{33} + a_{13}a_{21}a_{33} + a_{11}a_{22}a_{33} + a_{11}a_{23}a_{33} = 0, \\
a_{11}a_{22}a_{33} + a_{11}a_{23}a_{32} + a_{12}a_{22}a_{33} + a_{12}a_{23}a_{33} + a_{13}a_{22}a_{32} + a_{13}a_{22}a_{33} = 0, \\
a_{11}a_{22}a_{33} + a_{11}a_{22}a_{31} + a_{11}a_{21}a_{33} - a_{12}a_{21}a_{31} - a_{13}a_{21}a_{31} - a_{11}a_{23}a_{32} = 0, \\
a_{11}a_{22}a_{33} + a_{12}a_{23}a_{31} + a_{13}a_{21}a_{32} - a_{13}a_{22}a_{31} - a_{11}a_{23}a_{32} - a_{12}a_{21}a_{33} = 0, \\
a_{11}a_{22}a_{33} + a_{11}a_{22}a_{32} + a_{12}a_{22}a_{33} - a_{13}a_{22}a_{31} - a_{12}a_{21}a_{32} - a_{12}a_{23}a_{32} = 0, \\
a_{11}a_{22}a_{33} + a_{11}a_{23}a_{33} + a_{13}a_{22}a_{33} - a_{13}a_{23}a_{31} - a_{13}a_{23}a_{32} - a_{12}a_{21}a_{33} = 0. \\
\end{array}
\right.
$$
It~is straightforward to check that, apart from the $3\times 3$ universal matrices given in Example~\ref{ex:unversalMatricesInDim3},
the above equations are satisfied by the circulant matrix
\begin{equation}
\label{LogSpecificGood3x3-matrix}
C(0,1,-1)=
\left(
\begin{array}{rrr}
  0 &  1 & -1 \\
 -1 &  0 &  1 \\
  1 & -1 &  0 \\
\end{array}
\right)
\end{equation}
(and, by Lemma~\ref{lemma:permutationEquivMatrAreGood}, any of the permutation-similar matrices).
This matrix is specific to the choice of the logarithmic function.
It~cannot be obtained from the universal matrices by choosing suitable values of the parameters.
Any odd power of~(\ref{LogSpecificGood3x3-matrix}) is equal to~(\ref{LogSpecificGood3x3-matrix}) itself,
all of its even powers are equal and do not form good pairs with the logarithmic function.

We~remark that the inverse of the mapping defined by this matrix, that is, of the mapping
\begin{equation}
\label{logMapForCirculant(0,1,-1)}
\begin{array}{l}
f_1 = x_1 + \log(x_2-x_3),  \\
f_2 = x_2 + \log(x_3-x_1), \\
f_3 = x_3 + \log(x_1-x_2)   \\
\end{array}
\end{equation}
is not an elementary function and hence $\mathfrak{N}(f)=\infty.$ The change of variables $y_j = \exp(x_j),$ $g_j = \exp(\exp(f_j)),$ $j=1,2,3$
transforms~(\ref{logMapForCirculant(0,1,-1)}) into the equations
$$
g_1 = \left(\frac{y_2}{y_3}\right)^{y_1}, \quad g_2 = \left(\frac{y_3}{y_1}\right)^{y_2}, \quad g_3 = \left(\frac{y_1}{y_2}\right)^{y_3}.
$$

This example shows that the inverse of a Jacobian mapping in general need not be given by its finite Newton's inverse.
In~particular, the inverse of a Jacobian mapping of the form~(\ref{generalMapWithFcnPhi}) is in general
not a finite superposition of the function~$\varphi(\cdot)$ and arithmetic operations.
In~the present paper, (\ref{logMapForCirculant(0,1,-1)}) is the only example of a mapping with this property.
All the other mappings exposed here are inverted by a finite number of Newton's iterations~(\ref{NewtonIterations}).
\end{example}

Despite the analytic simplicity of the logarithmic function and the rationality of its derivative,
the structure of the set of all matrices which form good pairs with $\log(\cdot)$ is not clear at all.
For instance, the circulant matrices  $C(v,-v),$ for any $v\in\C^n,$
$C\left(0,1,-1-\sqrt{-1},\sqrt{-1}\right),$
$C\left(0,1+\sqrt{-1},-2,1-\sqrt{-1}\right),$
$C\left(0, 2, -1 - \sqrt{5}, 1 + \sqrt{5}, -2\right),$
and even
{\small
$$
\begin{array}{c}
C\left(0,\, \left(\sqrt{5}+1\right) \sqrt{10-2 \sqrt{5}}-2 \sqrt{-1} \left(\sqrt{5}+3\right),\,
\left(\sqrt{5}+3\right) \sqrt{10-2 \sqrt{5}}+4 \sqrt{-1} \left(\sqrt{5}+2\right), \phantom{--} \right. \\
\left.\phantom{------------}
-2 \left(\left(\sqrt{5}+2\right) \sqrt{10-2 \sqrt{5}}- \sqrt{-1} \left(\sqrt{5}+1\right)\right),\,
-4 \sqrt{-1} \left(\sqrt{5}+1\right)\right),
\end{array}
$$
}
as well as their complex conjugates form good pairs with the logarithmic function.
These circulant matrices are not particular instances of any universal matrices.

Another example of a family of matrices that form good pairs with $\log(\cdot)$ is given by
$$
\left(
\begin{array}{cccc}
 0 & a & -a-b & b \\
 a s & 0 & (-a-b) (s+t) & b s+a t+b t \\
 s & -s-t & 0 & t \\
 b s & -b s-a t-b t & (a+b) t & 0 \\
\end{array}
\right), \quad a,b,s,t\in\C.
$$
The problem of describing the set of all matrices which form good pairs with the logarithmic function
appears to have formidable complexity. We~remark that the property of forming a good pair with the logarithmic function
is in general not preserved by the Kronecker product as exemplified by the matrix
$
C(0,1,-1)
\otimes
\left(
\begin{array}{cc}
 1 & 1 \\
 1 & 0 \\
\end{array}
\right).
$

\begin{example}\rm
\label{ex:JMatrixWithCross-ShapedSupport}
For any $a,b\in(\C^{*})^{n-1}$ the matrix
$$
M=
\left(
\begin{array}{cccc}
0       & a_1    & \ldots & a_{n-1} \\
b_1     & 0      &      0 & 0       \\
\ldots  & \ldots & \ldots & \ldots  \\
b_{n-1} & 0      &      0 & 0       \\
\end{array}
\right)
$$
has rank~2 and forms a good pair with the function~$\zeta^d$ for $d=1,2,3,\ldots$ if and only if
\begin{equation}
\label{whenCrossSuppMatrixIsGood}
a_1 b_1^d + a_2 b_2^d + \ldots a_{n-1} b_{n-1}^d = 0.
\end{equation}
If~the condition~(\ref{whenCrossSuppMatrixIsGood}) is satisfied then the inverse of the mapping $f=x+(Mx)^d$ is given~by
$$
x=f-(M(f-(Mf)^d))^d.
$$
If~all of~$b_{j}$ are different then by Lemma~\ref{lemma:genSolToVandermondeSystem}
the matrix~$M$ cannot be universal no matter what $a\in(\C^{*})^{n-1}$ are.
If~some of the numbers~$b_{j}$ coincide, we may without loss of generality assume that
$$
b_1 = b_2 = \ldots = b_{j_1}, \,\,\, b_{j_1 + 1} = b_{j_1 + 2} = \ldots = b_{j_2}, \, \ldots, \,
b_{j_{k-1} + 1} = b_{j_{k-1} + 2} = \ldots = b_{n-1},
$$
and that no other elements in the set $\{b_1,\ldots,b_{n-1}\}$ are equal.
That is, we may assume that the elements of the first column of the matrix~$M$ are sorted to form groups consisting of equal elements.
This can be achieved by letting a suitable permutation act on the rows and columns of~$M.$
Under this assumption, the matrix~$M$ is universal if and only~if
$$
a_1 + a_2 + \ldots + a_{j_1} =  a_{j_1 + 1} + a_{j_1 + 2} + \ldots + a_{j_2} = \ldots = a_{j_{k-1} + 1} + a_{j_{k-1} + 2} + \ldots + a_{n-1} = 0.
$$
\end{example}

\begin{example}
\label{ex:goodMatr,dim=4,deg=3,rank=2,all homogeneities}\rm
For any $a,b,c,d,s_j\in\C$ the rank~$2$ matrix
$$
M=
\left(
\begin{array}{cccc}
 a\, s_1^3 s_2 s_3 s_4 &  b\, s_2^4 s_3 s_4     &  b\, s_2 s_3^4 s_4      &  a\, s_2 s_3 s_4^4     \\
 c\, s_1^4 s_3 s_4     &  d\, s_1 s_2^3 s_3 s_4 &  d\, s_1 s_3^4 s_4      &  c\, s_1 s_3 s_4^4     \\
-c\, s_1^4 s_2 s_4     & -d\, s_1 s_2^4 s_4     & -d\, s_1 s_2 s_3^3 s_4  & -c\, s_1 s_2 s_4^4     \\
-a\, s_1^4 s_2 s_3     & -b\, s_1 s_2^4 s_3     & -b\, s_1 s_2 s_3^4      & -a\, s_1 s_2 s_3 s_4^3 \\
\end{array}
\right)
$$
is permutation-similar to the Hadamard product
$$
H_{4,3}(s) \odot
\left(
\left(
\begin{array}{cc}
a & b \\
c & d \\
\end{array}
\right)
\otimes
\left(
\begin{array}{cc}
1 & -1 \\
1 & -1 \\
\end{array}
\right)
\right)
$$
and forms a good pair with the function~$\zeta^3$ i.e., for $x\in\C^4$ the mapping $x+(Mx)^3$ has unit Jacobian.
The inverse mapping is given by $f-(Mf)^3.$
\end{example}

\begin{example}
\label{ex:dim=4,deg=2,rank=2,rational parametrization} \rm
For any nondegenerate matrices $A=\left(\begin{array}{cc}a&b\\ c&d \end{array}\right)$ and $B=\left(\begin{array}{cc}s&t\\ u&v \end{array}\right)$
the matrix~$M$ of rational functions
\begin{equation}
\label{strangeRank2MatrixInDim=4}
{
\left(
\begin{array}{cccc}
 a       & b       & \frac{b d s u^2-a b s v u+b c t v u-a^2 t v^2}{(sv-tu) \left(b u s^2+d t u s+a t v s+c t^2 v\right)} & \frac{-b d u s^2+a b t u s-b c t v s+a^2 t^2 v}{(sv-tu) \left(b s u^2+d s v u+a t v u+c t v^2\right)} \\
 c       & d       & \frac{d^2 s u^2-b c s v u+c d t v u-a c t v^2}{(sv-tu) \left(b u s^2+d t u s+a t v s+c t^2 v\right)} & \frac{-d^2 u s^2+b c t u s-c d t v s+a c t^2 v}{(sv-tu) \left(b s u^2+d s v u+a t v u+c t v^2\right)} \\
 a s+c t & b s+d t & \frac{d u-a v}{sv-tu}                                                                                & \frac{t (as+ct)(atv+bsu) - s(bs+dt)(ctv+dsu)}{(sv-tu) \left(b s u^2+d s v u+a t v u+c t v^2\right)} \\
 a u+c v & b u+d v & \frac{u (bu+dv)(ctv+dsu)-v(au+cv)(atv+bsu)}{(sv-tu) \left(b u s^2+d t u s+a t v s+c t^2 v\right)}    & \frac{at-ds}{sv-tu} \\
\end{array}
\right)
}
\end{equation}
has rank~$2$ and forms a good pair with the function $\varphi(\zeta)=\zeta^2.$
Here we assume that all elements of~$M$ are well-defined, i.e., that none of the denominators vanish.
Of course, the denominators of the elements of~$M$ can be cleared by multiplying~$M$ with their least common multiple
but this will make~$M$ too cumbersome to display.

The first two rows of the matrix~(\ref{strangeRank2MatrixInDim=4}) are linearly independent for a generic choice of the matrix~$A$
while its last two rows are the linear combinations of the first two rows with the coefficients defined by the rows of the matrix~$B.$
The solution to the system of equations $f=x+(Mx)^2$ is given by the second iteration of the Newton mapping~(\ref{NewtonIterations}),
i.e., $x=f-(M(f-(Mf)^2))^2.$
\end{example}

\begin{example}
\label{ex:dim=8,deg=5,rank=2,two-dim family} \rm
For any $a,b\in\C$ the matrix
\begin{equation}
\label{dim=8,deg=5,rank=2,two-dim family}
{\tiny
M(a,b):=
\left(
\begin{array}{cccccccc}
 480 a & 0 & -360 a & -10 a & 32 a & 64 a & -9 a & -3 a \\
 0 & 480 b & -360 b & 10 b & 64 b & 32 b & -3 b & -9 b \\
 480 a & 480 b & -360 (a+b) & 10 (b-a) & 32 (a+2 b) & 32 (2 a+b) & 3 (-3 a-b) & 3 (-a-3 b) \\
 480 a & -480 b & -360 (a-b) & 10 (-a-b) & 32 (a-2 b) & 32 (2 a-b) & 3 (b-3 a) & 3 (3 b-a) \\
 480 a & 960 b & -360 (a+2 b) & 10 (2 b-a) & 32 (a+4 b) & 64 (a+b) & 3 (-3 a-2 b) & 3 (-a-6 b) \\
 960 a & 480 b & -360 (2 a+b) & 10 (b-2 a) & 64 (a+b) & 32 (4 a+b) & 3 (-6 a-b) & 3 (-2 a-3 b) \\
 1440 a & 480 b & -360 (3 a+b) & 10 (b-3 a) & 32 (3 a+2 b) & 32 (6 a+b) & 3 (-9 a-b) & -9 (a+b) \\
 480 a & 1440 b & -360 (a+3 b) & 10 (3 b-a) & 32 (a+6 b) & 32 (2 a+3 b) & -9 (a+b) & 3 (-a-9 b) \\
\end{array}
\right)
}
\end{equation}
has generic rank~$2$ and forms a good pair with the function $\varphi(\zeta)=\zeta^5.$
The normalized basis matrices $M(1,0)/240$ and $M(0,1)/240$ are idempotent.
The solution to the system of equations $f=x+(M(a,b)x)^5$ is given by the second iteration of the Newton mapping~(\ref{NewtonIterations}).
Similar examples can be constructed in any dimension $n\geq 4$ and for any degree of the monomial function~$\varphi(\zeta).$
\end{example}


\end{document}